\documentclass[10pt]{article}
\usepackage{amsmath,amssymb,amsthm,amscd,verbatim}

\begin{document}
\theoremstyle{plain}
\newtheorem{pro1}{Proposition}[section]
\newtheorem{cor1}[pro1]{Corollary}
\newtheorem{lem1}[pro1]{Lemma}
\newtheorem{thm1}[pro1]{Theorem}
\theoremstyle{definition}
\newtheorem{def1}[pro1]{Definition}
\newtheorem{fac1}[pro1]{Fact}
\newtheorem{rem1}[pro1]{Remark}
\newtheorem{ex1}[pro1]{Example}
\renewcommand{\labelenumi}{(\roman{enumi})}
\font\germ=eufm10
\def\al{\alpha}
\def\beneme{\begin{enumerate}}
\def\beq{\begin{equation}}
\def\beqn{\begin{eqnarray}}
\def\beqnn{\begin{eqnarray*}}
\def\bigsl{{\hbox{\fontD \char'54}}}
\def\cd{\cdots}
\def\del{\delta}
\def\Del{\Delta}
\def\ei{e_i}
\def\eit{\tilde{e}_i}
\def\eneme{\end{enumerate}}
\def\ep{\epsilon}
\def\eeq{\end{equation}}
\def\eeqn{\end{eqnarray}}
\def\eeqnn{\end{eqnarray*}}
\def\fit{\tilde{f}_i}
\def\ft{\tilde{f}}
\def\ge{{\mathfrak g}}
\def\gl{\hbox{\germ gl}}
\def\hom{{\hbox{Hom}}}
\def\ify{\infty}
\def\io{\iota}
\def\kp{k^{(+)}}
\def\km{k^{(-)}}
\def\llra{\relbar\joinrel\relbar\joinrel\relbar\joinrel\rightarrow}
\def\lan{\langle}
\def\lar{\longrightarrow}
\def\lm{\lambda}
\def\Lm{\Lambda}
\def\mapright#1{\smash{\mathop{\longrightarrow}\limits^{#1}}}
\def\nd{\noindent}
\def\nn{\nonumber}
\def\ot{\otimes}
\def\op{\oplus}
\def\opi{\ovl\pi_{\lm}}
\def\ovl{\overline}
\def\plm{\Psi^{(\lm)}_{\io}}
\def\qq{\qquad}
\def\q{\quad}
\def\qed{\hfill\framebox[3mm]{}}
\def\QQ{\mathbb Q}
\def\qi{q_i}
\def\qii{q_i^{-1}}
\def\ran{\rangle}
\def\rlm{r_{\lm}}
\def\ssl{\hbox{\germ sl}}
\def\slh{\widehat{\ssl_2}}
\def\ti{t_i}
\def\tii{t_i^{-1}}
\def\til{\tilde}
\def\tt{{\hbox{\germ{t}}}}
\def\ttt{\hbox{\germ t}}
\def\uq{U_q(\ge)}
\def\uqm{U^-_q(\ge)}
\def\uqp{U^+_q(\ge)}
\def\uqmq{{U^-_q(\ge)}_{\bf Q}}
\def\uqpm{U^{\pm}_q(\ge)}
\def\uqmp{U^{\mp}_q(\ge)}
\def\uqq{U_{\bf Q}^-(\ge)}
\def\uqz{U^-_{\bf Z}(\ge)}
\def\util{\tilde{U}_q(\ge)}
\def\vep{\varepsilon}
\def\vp{\varphi}
\def\vpi{\varphi^{-1}}
\def\xii{\xi^{(i)}}
\def\Xiioi{\Xi_{\io}^{(i)}}
\def\wtil{\widetilde}
\def\what{\widehat}
\def\wpi{\widehat\pi_{\lm}}
\def\ZZ{\mathbb Z}
\font\germ=eufm10
\def\ssl{\hbox{\germ sl}}
\def\slh{\widehat{\ssl_2}}
\title{\Large\bf Polyhedral Realizations of Crystal Bases 
\\ for Quantum Algebras of Finite Types}
\author{\normalsize AYUMU HOSHINO\thanks{e-mail address:
a-hoshin@mm.sophia.ac.jp \endgraf
{\q\negmedspace The author is supported
by Grant-in-Aid for JSPS Fellows.}} \\
\\
\normalsize Department of Mathematics, Sophia University, \\
\normalsize Tokyo 102-8554, JAPAN}
\date{}
\maketitle
\begin{abstract}
Polyhedral realization of crystal bases is one of the methods 
for describing the crystal base $B(\ify)$ explicitly.
This method can be applied to symmetrizable Kac-Moody types.
We can also apply this method to the crystal bases $B(\lm)$
of integrable highest weight modules and of modified quantum
algebras.
But, the explicit forms of the
polyhedral realizations of crystal bases $B(\ify)$ and $B(\lm)$
are only given in the case of arbitrary rank 2, of $A_n$ and of $A^{(1)}_n$.
So, we will give the polyhedral realizations of crystal bases $B(\ify)$
and $B(\lm)$
for all simple Lie algebras in this paper.
\end{abstract}
\tableofcontents

\renewcommand{\thesection}{\arabic{section}}
\section{Introduction}
\setcounter{equation}{0}
\renewcommand{\theequation}{\thesection.\arabic{equation}}

\q\, Quantum algebra $\uq:={\lan e_i, f_i, q^h \ran}_{i \in I}$
$(I= \{1,2,\cd,n\})$ which was introduced
in the study of solvable lattice models is applied
to the several kinds of
study of mathematical physics and
plays important roles.   
The nilpotent part $\uqm$$(={\lan f_i \ran}_{i \in I})$
of $\uq$ has a crystal base $B(\ify)$ which was constructed
by Kashiwara \cite{K1} 
and the irreducible integrable highest weight representaion
of $\uq$ also has crystal base $B(\lm)$.


The crystal base has been realized by several methods but it is not so easy
to obtain the explicit form.
Polyhedral realization of crystal bases is one of the methods for
realizing crystal base explicitly, which was introduced
by Nakashima and Zelevinsky \cite{NZ}.
We can describe a vector of crystal base $B(\ify)$ as a lattice point of 
certain convex polyhedron in an infinite $\mathbb Z$-lattice 
by this method. This method can be applied to
symmetrizable Kac-Moody types and applied to the crystal base $B(\lm)$
of the irreducible integrable highest weight module. 
In \cite{NZ}, polyhedral realizations of $B(\ify)$ are given when
$\ge$ is of arbitrary rank 2 cases, of $A_n$ and of $A^{(1)}_{n-1}$ and
in \cite{N1}, Nakashima gave the polyhedral realizations of the crystal base 
$B(\lm)$ $(\lm \in P_+)$ of irreducible 
integrable highest weight modules when $\ge$ is the same cases as the above. 
He and the author \cite{HN} applied this method 
to the modified quantum algebras
and had the polyhedral realizations of $B(\uq a_{\lm})$ and 
the some specific
connected component of $B_0(\lm)$ $(\lm \in P)$ containing 
$u_{\ify} \ot t_{\lm} \ot u_{-{\ify}}$ for $\ge$ is of type $A_n$
under certain assumption
on the weight $\lm$ and for $\ge$ is of type $A^{(1)}_1$ 
on positive level $\lm$ and
in \cite{H}, the author gave the one of $B_0(\lm)$
for arbitrary rank 2 cases on positive level $\lm$.
After Nakashima and Zelevinsky's work, P. Littelmann described
crystal base $B(\ify)$ by some inequalities
(which are called ``cone'') for all simple Lie algebras
and $B(\lm)$ for classical Lie algebras \cite{Li}.
In this paper, we will give the polyhedral realizations of
crystal bases $B(\ify)$ and $B(\lm)$ for all simple Lie algebras. 
In order to
treat the cases, we improve the Theorem in \cite{NZ} and
obtain the polyhedral realizations.

\vskip5pt
This paper is organized as follows:
in Section 2, we review the theory of crystal base and
methods of polyhedral realizations of $B(\ify)$ and $B(\lm)$.
In Section 3, we improve the method which is given by \cite{NZ}
and obtain the explicit forms of the polyhedral realizations
of $B(\ify)$ and $B(\lm)$ for all simple Lie algebras.
But, we do not write in this paper about the $B(\lm)$ for
$E_7$ and $E_8$
since numerous inequalities appear.
\vskip5pt

{\bf Acknowledgements.} The author would like to thank
Professor Toshiki Nakashima for his support and helpful comments.

\section{Preliminaries}
\setcounter{equation}{0}
\renewcommand{\thesection}{\arabic{section}}
\renewcommand{\theequation}{\thesection.\arabic{equation}}

\subsection{Crystal bases and Crystals}

In this subsection, we review the theory of
crystal bases and crystals. 
We fix a finite index set $I$ and
let $A=(a_{ij})_{i,j\in I}$ be 
a generalized symmetrizable Cartan matrix,
$(\ttt,\{\al_i\}_{i\in I},\{h_i\}_{i\in I})$ be the associated 
Cartan data and $\ge$ be the associated Kac-Moody Lie algebra
where $\al_i$ (resp. $h_i$) is called a simple root
(resp. simple coroot).
Let $P$ be a weight lattice with a $\QQ$-valued symmetric 
bilinear form $(\; , \;)$, $P^*$ be a dual lattice 
including $\{h_i\}_{\in I}$ and 
$Q:=\bigoplus_{i\in I}\QQ(q)\al_i$ be a root lattice.
The quantum algebra $\uq$ to be the associative algebra with 
$1$ over $\QQ(q)$ generated by $e_i,f_i, q^h$ ($i\in I,h\in P^*$) 
with the usual relations.
%
%
%
%
%
%
Let $\uqm := \lan f_i \ran_{i \in I}$ be the subalgebra
of $\uq$ and $V(\lm)$ be the irreducible integrable highest
weight module. $\uqm$ (resp.$V(\lm)$) have a
crystal base $(L(\ify),B(\ify))$ (resp.$(L(\lm),B(\lm))$)
satisfying some properties.
Let $\pi_{\lm} : \uqm \rightarrow V(\lm) \cong \uqm / 
\sum_i\uqm\fit^{1+\lan h_i, \lm \ran}$ be the canonical
projiction and $\hat{\pi}_{\lm} : L(\ify)/qL(\ify)
\rightarrow L(\lm)/qL(\lm)$ be the induced map from $\pi_{\lm}$.
We note that $\hat{\pi}_{\lm}(B(\ify)) = B(\lm) \sqcup \{0\}$.

The notion of crystal is obtained by
abstracting the combinatorial properties of crystal bases.
A crystal $B$ has maps $wt : B \rightarrow P$,
$\vep_i,\,\vp_i : B \rightarrow \ZZ \sqcup \{-\ify\}$ and
$\eit,\,\fit : B \sqcup \{0\}\rightarrow \ B \sqcup \{0\}$
with some aximos. In fact, crystal bases
$B(\ify)$ and $B(\lm)$ are also crystal.
The tensor product of the crystal bases is again crystal base
and so, we can consider the tensor product of crystals.

%


\subsection{Polyhedral Realization of $B(\ify)$}

In this subsection, 
we review the polyhedral realization of the crystal $B(\ify)$ 
(see \cite{NZ}).

First, we recall the crystal structure of $\ZZ^{\ify}$.
We consider the following additive groups:
\begin{eqnarray*}
\ZZ^{\ify}&:=&\{(\cd,x_k,\cd,x_2,x_1)\,|\, x_k\in\ZZ
\,\,{\rm and}\,\,x_k=0\,\,{\rm for}\,\,k\gg 0\}.
\end{eqnarray*}
We will denote by
 $\ZZ^{{\ify}}_{\geq 0} \subset \ZZ^{{\ify}}$ 
the semigroup of nonnegative sequences. 
Take an infinite sequence of indices
$\io=(\cd,i_k,\cd,i_2,i_1)$ 
from $I$ such that
\begin{equation}
{\hbox{
$i_k\ne i_{k+1}$ for any $k$, \,and $\sharp\{k>0
\,:\, i_k=i\}=\ify$ for any $i\in I$.}}
\label{seq-con}
\end{equation}
The crystal structure on $\ZZ^{\ify}$ associate to 
$\io$ is defined as follows. Let
$\vec x = (\cd,x_2,x_1) \in \ZZ^{\ify}$. We set for $k \geq 1$
\begin{eqnarray}
\sigma_k(\vec x) := x_k + \sum_{j >k}\lan h_{i_k}, \alpha_{i_j}\ran
x_j.
\end{eqnarray}
Since $x_j=0$ for $j\gg0$,
 $\sigma_k$ is
well-defined.
Let $\sigma^{(i)} (\vec x)
 := {\rm max}_{k: i_k = i}\sigma_k (\vec x)$ and
\begin{equation}
M^{(i)} = M^{(i)} (\vec x) :=
\{k: i_k = i, \sigma_k (\vec x) = \sigma^{(i)}(\vec x)\}.
\label{m(i)}
\end{equation}
Note that 
$\sigma^{(i)} (\vec x)\geq 0$, and that 
$M^{(i)} = M^{(i)} (\vec x)$ is finite set if and only if 
$\sigma^{(i)} (\vec x) > 0$. 
Now, we define the map 
$\eit: \ZZ^{\ify} \lar \ZZ^{\ify} \sqcup\{0\}$
, 
$\fit: \ZZ^{\ify} \lar \ZZ^{\ify}$, by
 $\eit(0)=\fit(0)=0$ and 
\begin{eqnarray}
(\fit(\vec x))_k &=& x_k + \delta_{k,{\rm min}\,M^{(i)}},
\label{action-f}\\
(\eit(\vec x))_k &=& 
x_k - \delta_{k,{\rm max}\,M^{(i)}} \,\, {\rm if}\,\,
\sigma^{(i)}(\vec x)>0; \;\,
 {\text{ other wise}} \,\, \eit(\vec x)=0,
\label{action-e}
\end{eqnarray}
where $\del_{i,j}$ is Kronecker's delta.
We also define the weight function and the function
$\vep_i$ and $\vp_i$ on $\ZZ^{\ify}$ as follows:
\begin{equation}
\begin{array}{l}
wt(\vec x) := -\sum_{j=-{\ify}}^{\ify} x_j \al_{i_j}, \,\,
\vep_i (\vec x) := \sigma^{(i)} (\vec x), \,\,
\vp_i (\vec x) := \lan h_i, wt(\vec x) \ran + \vep_i(\vec x).
\end{array}
\label{wt-vep-vp}
\end{equation}
We denote this crystal by $\ZZ^{\ify}_{\io}$.


\begin{pro1}[\cite{K2}]
There is a unique embedding of crystals according to $\io$
\begin{eqnarray}
\Psi_{\io}:&B(\ify)& \hookrightarrow \; \ZZ^{{\ify}}_{\geq 0}
\subset \ZZ^{{\ify}}_{\io}
\q
(u_{{\ify}}\mapsto (\cd,0,0)).
\label{kas+}
\end{eqnarray}
\end{pro1}
\noindent
This embedding is called {\it Kashiwara embedding}.

Next, we review the polyhedral realization of $B(\ify)$ for
describing the image of {\it Kashiwara embedding}.
We consider the following infinite dimensional vector spaces
and their dual spaces:
\begin{eqnarray*}
\QQ^{{\ify}}&:=&\{\vec{x}=
(\cd,x_k,\cd,x_2,x_1): x_k \in \QQ\,\,{\rm and }\,\,
x_k = 0\,\,{\rm for}\,\, k \gg 0\},\\
 (\QQ^{{\ify}})^*&:=&{\rm Hom}(\QQ^{{\ify}},\QQ).
\end{eqnarray*}
We will write a linear form $\vp \in (\QQ^{\ify})^*$ as
$\vp(\vec{x})=\sum_{k \geq 1} \vp_k x_k$ ($\vp_j\in \QQ$). 
For the sequence 
$\io=(i_k)_{k \geq 1}$ and 
$k \geq 1$, we set 
\[
\kp:={\rm min}\{l:l>k>0\,\,\,\,
{\rm  and }\,\,i_k=i_l\},
\]
if it exists, otherwise $\km=0$.
We define a linear form $\beta_k$ $(k\geq1)$ on $\QQ^{{\ify}}$ by
\begin{eqnarray}
\beta_{k}(\vec x)
&=& \sigma_k(\vec x) - \sigma_{k^{(+)}}(\vec x) \\
&=&
 \, x_k+\sum_{k<j<\kp}\lan h_{i_k},\al_{i_j}\ran x_j+x_{\kp}. \nn
\label{beta}
\end{eqnarray}
By using these linear forms, let us
 define a piecewise-linear operator 
$S_k=S_{k,\io}$ on $(\QQ^{{\ify}})^*$
as follows:
\begin{equation}
S_k(\vp):=\begin{cases}
\vp-\vp_k\beta_k & \text{if}\;\, \vp_k>0,\\
 \vp-\vp_k\beta_{\km} & \text{if}\;\, \vp_k\leq 0,\\
\end{cases}
\label{Sk}
\end{equation}
for $\vp(\vec x)=\sum \vp_k x_k\in (\QQ^{{\ify}})^*$.
Here we set
\begin{eqnarray}
\Xi_{\io} &:= &
\{S_{j_l}\cd S_{j_2}S_{j_1}( x_{{j_0}})\,|\,
l\geq0,j_0,j_1,\cd,j_l\geq1\},
\label{xiiota}\\
\Sigma_{\io} & := &
\{\vec x\in \ZZ^{{\ify}}\subset \QQ^{{\ify}}
\,|\,\vp(\vec x)\geq0\,\,{\rm for}\,\,
{\rm any}\,\,\vp\in \Xi_{\io}\}.
\end{eqnarray}
We impose on $\io$ the following ``positivity assumption'' (P):
\begin{eqnarray}
&{\hbox{(P)\;   for $\io$, if $\km=0$ then $\vp_k\geq0$ for any 
$\vp(\vec x)=\sum_k\vp_kx_k\in \Xi_{\io}$}}.\label{p}
\end{eqnarray}

\begin{thm1}[\cite{NZ}]
Let $\io$ be the sequence of indices satisfying 
$(\ref{seq-con})$ and the positivity assumption $(P)$. 
Let $\Psi_{\io}:B(\ify)\hookrightarrow \ZZ^{\ify}_{\io}$
be the 
Kashiwara embedding. Then, we have
${\rm Im}(\Psi_{\io})(\cong B(\ify))=\Sigma_{\io}$.
\label{nz}
\end{thm1}
\noindent
We call $\Sigma_{\io}$ the polyhedral
realization of $B(\ify)$.

\vskip5pt

\subsection{Polyhedral Realization of $B(\lm)$}

In this subsection, 
we review the polyhedral realization of the crystal $B(\lm)$ 
(see \cite{N1}). 

Let $R_{\lm}:= \{r_{\lm} \}$ be the crystal for $\lm \in P_+$
($P_+$: set of dominant integral weights)
defined in \cite{N1}. For the crystal $B(\ify) \ot R_{\lm}$,
we define the following map:
\begin{eqnarray}
\Phi_{\lm} : (B(\ify) \ot R_{\lm}) \sqcup \{0\} \rightarrow B(\lm)
\sqcup \{0\}
\end{eqnarray}
by $\Phi_{\lm}(0) = 0$ and $\Phi_{\lm}(b \ot r_{\lm}) = \hat\pi_{\lm}(b)$
for $b \in B(\ify)$. We set
\begin{eqnarray*}
\tilde{B}(\lm) := \{b \ot r_{\lm} \in B(\ify) \ot R_{\lm} : 
\Phi_{\lm}(b \ot r_{\lm}) \ne 0 \}.
\end{eqnarray*}
\begin{thm1}[\cite{N1}]
(i) The map $\Phi_{\lm}$ becomes a surjective strict morphism
of crystals $B(\lm) \ot R_{\lm} \rightarrow B(\lm)$.

(ii) $\tilde{B}(\lm)$ is a subcrystal of $B(\lm) \ot R_{\lm}$,
and $\Phi_{\lm}$ induces the isomorphism of crystals
$\tilde{B}(\lm) \simeq B(\lm)$.
\label{nint}
\end{thm1}
 
Let us denote $\ZZ_{\io} \ot R_{\lm}$ by $\ZZ_{\io}[\lm]$.
Here note that we can identify $\ZZ_{\io}[\lm]$ with
$\ZZ_{\io}^{\ify}$ as a set since the crystal $R_{\lm}$ has 
only one element but their crystal structures are different.
By Theorem \ref{nint}, we have the strict embedding of crystals
$\Omega_{\lm} : B(\lm) (\cong \tilde{B}(\lm)) \hookrightarrow 
B(\ify) \ot R_{\lm}$. Combining $\Omega_{\lm}$ and Kashiwara
embedding $\Psi_{\io}$, we obtain the following:

\begin{thm1}[\cite{N1}]
There exists the unique strict embedding of crystals
\begin{eqnarray*}
\Psi_{\io}^{(\lm)} : B(\lm) \hookrightarrow B(\ify) \ot R_{\lm}
\hookrightarrow \ZZ_{\io}^{\ify} \ot R_{\lm} =: \ZZ_{\io}^{\ify}[\lm]
\end{eqnarray*}
such that $\Psi_{\io}^{(\lm)}(u_{\lm}) = (\cd,0,0) \ot r_{\lm}$.
\label{embed}
\end{thm1}
We fix an infinite sequence of indices
$\io=(\cd,i_k,\cd,i_2,i_1)$ 
satisfying (\ref{seq-con}) and $\lm \in P_+$.
We define a linear form $\beta_k^{(\pm)}$ $(k\geq1)$ on $\QQ^{{\ify}}$ by
\begin{align}
\beta_{k}^{(+)}(\vec x)
\,&= \, x_k+\sum_{k<j<\kp}\lan h_{i_k},\al_{i_j}\ran x_j+x_{\kp},\\
\beta_{k}^{(-)}(\vec x)
\,&=
\begin{cases}
x_{k^{(-)}}+\sum_{k^{(-)}<j<k}\lan h_{i_k},\al_{i_j}\ran x_j+x_k,
& \text{ if }\; k^{(-)} >0, \\
- \lan h_{i_k}, \lm \ran +
\sum_{1<j<k}\lan h_{i_k},\al_{i_j}\ran x_j+x_k,
& \text{ if }\; k^{(-)} =0.
\end{cases}
\label{betaint}
\end{align}
Here note that $\beta_k^{(+)} = \beta_k$ and $\beta_k^{(-)} = \beta_{k^{(-)}}$
if $k^{(-)} >0$.
By using these linear forms, let us
 define a piecewise-linear operator 
$\hat{S}_k=\hat{S}_{k,\io}$ on $(\QQ^{{\ify}})^*$
as follows:
\begin{equation}
\hat{S}_k(\vp):=\begin{cases}
\vp-\vp_k\beta_k^{(+)} & \text{if}\;\, \vp_k>0,\\
 \vp-\vp_k\beta_k^{(-)} & \text{if}\;\, \vp_k\leq 0,\\
\end{cases}
\label{Skhat}
\end{equation}
for $\vp(\vec x)=c+\sum \vp_k x_k$
$(c,\vp_k \in \QQ)$ on $\QQ^{{\ify}}$.
For the fixed sequence $\io = (i_k)$, in case $k^{(-)}=0$ for
$k \geq 1$, there exists unique $i \in I$ such that $i_k = i$.
We denote such $k$ by $\io^{(i)}$, namely, $\io^{(i)}$ is the
first number $k$ such that $i_k = i$. Here we set for $\lm \in P_+$
and $i \in I$
\begin{eqnarray}
\lm^{(i)}(\vec x) := -\beta_{\io^{(i)}}^{(-)}(\vec x)
= \lan h_i, \lm \ran - 
\sum_{1<j<\io^{(i)}}\lan h_i,\al_{i_j}\ran x_j-x_{\io^{(i)}}.
\end{eqnarray}
For $\io$ and $\lm \in P_+$, let $\Xi_{\io}[\lm]$
be the set of all linear functions generated by $\hat{S}_k$
from the functions $x_j$ $(j\geq 1)$ and $\lm^{(i)}$ $(i \in I)$,
namely,
\begin{eqnarray*}
\Xi_{\io}[\lm] &:= &
\{\hat{S}_{j_l}\cd \hat{S}_{j_2}\hat{S}_{j_1}( x_{{j_0}}) :
l\geq0,j_0,j_1,\cd,j_l\geq1\} \\
&&\cup \{\hat{S}_{j_k}\cd \hat{S}_{j_2}\hat{S}_{j_1}(\lm^{(i)}(x)) :
k\geq0,\; i \in I, \;j_0,j_1,\cd,j_k\geq1\}.
\end{eqnarray*}
Now we set
\begin{eqnarray}
\Sigma_{\io} & := &
\{\vec x\in \ZZ_{\io}^{{\ify}}[\lm]\subset \QQ^{{\ify}}
: \vp(\vec x)\geq0\,\,{\rm for}\,\,
{\rm any}\,\,\vp\in \Xi_{\io}[\lm]\}.
\end{eqnarray}
For $\io$ and $\lm \in P_+$, a pair $(\io, \lm)$
is called $ample$ if $\Sigma_{\io}[\lm] \ni \vec 0 = (\cd,0,0)$. 

\begin{thm1}[\cite{N1}]
Suppose that $(\io, \lm)$ is ample. Let
$\Psi_{\io}^{\lm} : B(\lm) \hookrightarrow \ZZ_{\io}^{\ify}[\lm]$
be the embedding as in Theorem \ref{embed}. Then we have
${\rm Im}(\Psi_{\io}^{(\lm)})(\cong B(\lm))=\Sigma_{\io}[\lm]$.
\label{n1}
\end{thm1}
\noindent
We call $\Sigma_{\io}[\lm]$ the polyhedral
realization of $B(\lm)$.

In the rest of this subsection, we note about $\hat{S}_k$ and $S_k$.
We define the linear form $\xi^{(i)}$ $(i \in I)$ on $\QQ$ by
\begin{eqnarray}
\xi^{(i)}(\vec x) :=
- \sum_{1<j<\io^{(i)}}\lan h_i,\al_{i_j}\ran x_j-x_{\io^{(i)}}
= - \lan h_i, \lm \ran + \lm^{(i)}(\vec x)
\end{eqnarray}
and set of linear forms $\Xi_{\io}^{(i)}$ by
\begin{eqnarray}
\Xi_{\io}^{(i)} := \{S_{j_l}\cd S_{j_2}S_{j_1}\xi^{(i)} :
l\geq0,j_0,j_1,\cd,j_l\geq1\}.
\end{eqnarray}
Here we introduce the {\it strict positivity assumption} for
$\io$ as follows:
\begin{eqnarray}
&{\text{if $\km=0$ then $\vp_k\geq0$ for any 
$\vp(\vec x)=\sum_k\vp_kx_k\in (\Xi_{\io}$}\bigcup \cup_{j \in I}
\Xi_{\io}^{(j)})\setminus \{\xi^{(i)} : i \in I \} },\label{sp}
\end{eqnarray}
where $\Xi_{\io}$ is defined by (\ref{xiiota}). Then we have
following Lemma:
\begin{lem1}[\cite{N1}]
Under the strict positivity assumption for $\io$, we have
\begin{eqnarray}
\hat{S}_{j_l}\cd \hat{S}_{j_2}\hat{S}_{j_1}x_{j_0} 
= S_{j_l}\cd S_{j_2}S_{j_1}x_{j_0},
\end{eqnarray}
for any $l \geq 0$, $j_0, j_1, \cd, j_l \geq 1$, and
\begin{eqnarray}
\hat{S}_{j_l}\cd \hat{S}_{j_2}\hat{S}_{j_1}\lm^{(i)}(\vec x) 
= \lan h_i, \lm \ran + S_{j_l}\cd S_{j_2}S_{j_1}\xi^{(i)}(\vec x),
\label{abcd}
\end{eqnarray}
for any $l \geq 0$, $j_0, j_1, \cd, j_l \geq 1$ and $i \in I$,
if the left hand side of (\ref{abcd}) is non zero.
\label{lmxi}
\end{lem1}

\section{Explicit forms of polyhedral realization of 
$B(\ify)$ and $B(\lm)$}
\setcounter{equation}{0}
\renewcommand{\theequation}{\thesection.\arabic{equation}}

The polyhedral realizations of $B(\ify)$ and $B(\lm)$ are already known
in the case of all rank $2$ Kac-Moody types, 
of $A_n$ and of $A^{(1)}_{n-1}$ (\cite{NZ}, \cite{N1}). 
In this section, we treat other simple types.
We obtain the following Theorem by \cite{NZ}
Theorem 3.1:

\begin{thm1}
Let $\Xi_{\io}'$ be a set of the linear forms in $\QQ^{\ify}$ and set
\begin{eqnarray*}
\Sigma_{\io}' := \{ \vec x \in \ZZ^{\ify}_{\io}\;|\;\vp(x)\geq 0
\text{ for any } \vp \in \Xi_{\io}' \}
\end{eqnarray*} 
satisfying the following conditions:

\begin{enumerate}
\item $\Xi_{\io}'$ is closed under the action of $S_k$'s,
\item $\io$ satisfies the positivity assumption (P),
\item all entries of $\vec x = (\cd,x_2,x_1) \in  
\Sigma_{\io}'$ are nonnegative.
\end{enumerate}
Then we have
\begin{eqnarray*}
{\rm Im}(\Psi_{\io})(\cong B(\ify))=\Sigma_{\io}'.
\end{eqnarray*}
\begin{proof}
We recall that the following two facts for certain $\Sigma_{\io}
(\subset \ZZ_{\io}^{\ify})$  (see \cite{NZ} Theorem 3.1):

\noindent
$A)$ If $\Sigma_{\io}$ is closed under the action of $\fit$,
then Im $(\Psi_{\io}) \subset \Sigma_{\io}$,

\noindent
$B)$ If $\Sigma_{\io}$ is closed under the action of $\eit$
and all the entries of $\vec x \in \Sigma_{\io}$ are nonnegative,
then $\Sigma_{\io} \subset \text{Im} (\Psi_{\io})$ since
every $\vec x \in \Sigma_{\io}$ can be transformed to $\vec 0$
by the action of $\eit$'s ($i \in I$).

We assume the condition $(i)$. The condition $(i)$ means that
$\Sigma_{\io}'$ is closed under the action of $\eit$'s.
Using $(ii)$, we obtain that 
$\Sigma_{\io}'$ is closed under the action of $\fit$'s
and this shows Im$(\Psi_{\io}) \subset \Sigma_{\io}'$.
Using $(iii)$, we obtain that 
$\Sigma_{\io}' \subset \text{Im} (\Psi_{\io})$.

\end{proof}
\label{nz_modi}
\end{thm1}

\noindent
\begin{rem1}Theorem \ref{nz_modi} implies that we need
not consider a set of linear functions $\Xi_{\io}$ which is closed under the
action of $S_k$ to $x_{j}$ for any $j>0$. We only need to find
a set of linear functions $\Xi_{\io}$ which is closed under the action of
$S_k$ to $x_{j}$ for some $j>0$ and which satisfies two conditions
$(ii)$, $(iii)$.
\end{rem1}

It will be convenient for us to change the indexing set for
$\ZZ^{\ify}$ from $\ZZ_{\geq 1}$ to 
$\ZZ_{\geq 1} \times [1,n]$. We will do this with the help of the 
bijection $\ZZ_{\geq1} \times [1,n] \to \ZZ_{\geq 1}$ given by  
($(j;i) \mapsto (j-1)n + i$). Thus, we will write an element
 $\vec x \in \ZZ^{\ify}$ 
as doubly-indexed family $(x_{j;i})_{j \geq 1, i \in [1,n]}$ of
nonnegative integers.
Therefore, we can write that $\vec x \in \ZZ^{\ify}$ as   
$(\cd,x_{1;2},x_{1;1})$.
We will adopt the convention that $x_{j;0} = x_{j;n+1} = 0$
unless $i \in [1,n]$. Hereafter, we fix the infinite sequence $\io$
as follows:
$$\io := (\cd,\underbrace{n,n-1,\cd,2,1},\cd,
\underbrace{n,n-1,\cd,2,1},\underbrace{n,n-1,\cd,2,1}),$$
where $n$ is the size of Cartan matrix.

\subsection{$B_n$ case}

We consider the case of type $B_n$ and give the
explicit forms of $\Sigma_{\io}$ and $\Sigma_{\io}[\lm]$.
First, we give the polyhedral realization $\Sigma_{\io}$.
We define for any $j \geq 1,\, 0 \leq k \leq
2n-1$ 
\begin{eqnarray}
\vp_{j;k}:=
\begin{cases}
  id & (k=0),\\
  S_{j;k}S_{j;k-1} \cd S_{j;2}S_{j;1} & (1 \leq k \leq n-1),\\
  S_{j+k-n;2n-k}S_{j+k-1-n;2n-k+1} \cd S_{j+1;n-1}S_{j;n}\vp_{j;n-1}
 & (n \leq k \leq 2n-1).
\end{cases}
\label{c_phi_s}
\end{eqnarray}

\begin{lem1}
\begin{eqnarray}
\vp_{j;k}(x_{j;1})=
\begin{cases}
x_{j;k+1}-x_{j+1;k} & (0 \leq k \leq n-1),\\
x_{j+k-n+1;2n-k-1}-x_{j+k-n+1;2n-k} & (n \leq k \leq 2n-1).
\end{cases}
\label{c_phi_x}
\end{eqnarray}

\begin{proof}
By the induction on $k$. The case of $k=0$ is trivial.\\

\noindent
Case I) $1 \leq k \leq n-1$.\\

\noindent
If $k=1$, we have
\begin{eqnarray*}
S_{j;1}(x_{j;1}) &=& x_{j;1} - (x_{j;1}-x_{j;2}+x_{j+1;1}) \\
		 &=& x_{j;2}-x_{j+1;1} \\
		 &=& \vp_{j;1}(x_{j;1}).
\end{eqnarray*}

\noindent
If $1 < k \leq n-1$, we assume that 
$\vp_{j;k-1}(x_{j;1}) = x_{j;k}-x_{j+1;k-1}$. Then we have
\begin{eqnarray*}
S_{j;k}\vp_{j;k-1}(x_{j;1}) &=& S_{j;k}(x_{j;k}-x_{j+1;k-1}) \\ 
        &=& x_{j;k}-x_{j+1;k-1}-(x_{j;k}-x_{j;k+1}-x_{j+1;k-1}+x_{j+1;k}) \\
	&=& x_{j;k+1}-x_{j+1;k} \\
	&=& \vp_{j;k}(x_{j;1}).
\end{eqnarray*}

\noindent
Case II) $n \leq k \leq 2n-1$.\\

\noindent
If $k=n$, using the result of Case I), we have
\begin{eqnarray*}
S_{j;n}\vp_{j;n-1}(x_{j;1}) &=& S_{j;n}(x_{j;n}-x_{j+1;n-1}) \\
	&=& x_{j;n}-x_{j+1;n-1}-(x_{j;n}-2x_{j+1;n-1}+x_{j+1;n}) \\
	&=& x_{j+1;n-1}-x_{j+1;n} \\
	&=& \vp_{j;n}(x_{j;1}).
\end{eqnarray*}

\noindent
If $n < k \leq 2n-1$, we assume that 
$\vp_{j;k-1}(x_{j;1}) = x_{j+k-n;2n-k}-x_{j+k-n;2n-k+1}$. Then we have
\begin{eqnarray*}
S_{j+k-n;2n-k}\vp_{j;k-1}(x_{j;1}) &=& S_{j+k-n;2n-k}
(x_{j+k-n;2n-k}-x_{j+k-n;2n-k+1}) \\ 
&=& x_{j+k-n;2n-k}-x_{j+k-n;2n-k+1} \\
&&-(x_{j+k-n;2n-k}-x_{j+k-n;2n-k+1}-x_{j+k-n+1;2n-k-1}+x_{j+k-n+1;2n-k}) \\
&=& x_{j+k-n+1;2n-k-1}-x_{j+k-n+1;2n-k} \\
&=& \vp_{j;k}(x_{j;1}).
\end{eqnarray*}
\end{proof}
\label{c1}
\end{lem1}

\begin{lem1}
$\vp_{j;k}(x_{j;1})$ is closed
under the actions of all transformations $S_{m;l}$ for
any $m \geq 1, l \geq 1$.
\begin{proof}
Using the definition of $\vp_{j;k}$ (\ref{c_phi_s}) and the formula
(\ref{c_phi_x}) in Lemma \ref{c1},
if $k=0$, $\vp_{j;0} = id$ and we have
\begin{eqnarray*}
S_{m;l}(\vp_{j;0}(x_{j;1}))=
\begin{cases}
\vp_{j;1}(x_{j;1}) & if \;\;  (m;l)=(j;1),\\
x_{j;1} & other \; wise.
\end{cases}
\end{eqnarray*}

If $1 \leq k \leq n-1$, we have
\begin{eqnarray*}
S_{m;l}(\vp_{j;k}(x_{j;1}))=
\begin{cases}
\vp_{j;k-1}(x_{j;1}) & if \;\;  (m;l)=(j;k),\\
\vp_{j;k+1}(x_{j;1}) & if \;\;  (m;l)=(j;k+1),\\
\vp_{j;k}(x_{j;1}) & other \; wise.
\end{cases}
\end{eqnarray*}

If $n \leq k \leq 2n-1$, we have
\begin{eqnarray*}
S_{m;l}(\vp_{j;k}(x_{j;1}))=
\begin{cases}
\vp_{j;k-1}(x_{j;1}) & if \;\;  (m;l)=(j+k-n;2n-k),\\
\vp_{j;k+1}(x_{j;1}) & if \;\;  (m;l)=(j+k-n+1;2n-k-1),\\
\vp_{j;k}(x_{j;1}) & other \; wise.
\end{cases}
\end{eqnarray*}

\end{proof}
\label{c2}
\end{lem1}

Now, we define
\begin{eqnarray*}
\Xi_{\io} &:=& \{\; \vp_{j;k}(x_{j;1}) :
j \geq 1, \; 0 \leq k \leq 2n-1 \}, \\
\Sigma_{\io} &:=& \{\; \vec x \in \ZZ^{\ify}_{\io} :
\vp(\vec x) \geq 0 \text{ for any } \vp \in \Xi_{\io} \;\}.
\end{eqnarray*}

\begin{thm1}
\label{c3}
Polyhedral realization $\Sigma_{\io}$
of $B(\ify)$ for type $B_n$ is described as
below:
\begin{eqnarray*}
&&x_{j;i} = 0 \q \text{ for }\; j, i \notin [1,n], \\
&&x_{1;i}\geq x_{2;i-1}\geq \cd \geq x_{i;1} \geq 0
 \q \text{ for }\; 1 \leq i \leq n-1,\\
&&x_{j;n}\geq x_{j+1;n-1}\geq \cd \geq x_{n;j}\geq 0
\q \text{ for }\; 1 \leq j \leq n,\\
&&x_{j;n-j+1} \geq x_{j;n-j+2} \geq \cd \geq x_{j;n} \geq 0
\q \text{ for }\; 2 \leq j \leq n.
\end{eqnarray*}

\begin{proof}
First, we show that $\Sigma_{\io}$ is the polyhedral realization of
$B(\ify)$,
so we shall check the conditions of Theorem \ref{nz_modi}.
$\Xi_{\io}$ is closed under the action of $S_k$ by Lemma \ref{c2}.
The coefficients of $x_{1;i}$ $(i = 1,2,\cd,n)$ are
positive for $\vp \in \Xi_{\io}$ by Lemma $\ref{c1}$, and
this shows that $\io$ satisfies the ``positivity assumption''.


We will show that all entries of $\vec x \in  \Sigma_{\io}$
are nonnegative.
In the case of $m \geq 1$, $0 \leq l \leq n-1$ for
$\vp_{m;l}(x_{m;1})$, we have 
$$ x_{m;l+1} \geq x_{m+1;l}$$
and we consider the cases that 
$(m;l)=(j;i-1),(j+1;i-2),\cd,(j+i-2;1),(j+i-1;0)$
for any $j \geq 1$, $1 \leq i \leq n$, then we obtain
$$x_{j;i}\geq x_{j+1;i-1}\geq \cd
\geq x_{j+i-2;2}\geq x_{j+i-1;1} \geq 0.$$
This shows that $x_{j;i} \geq 0$ for any $j \geq 1$, $1 \leq i \leq n$.
Therefore, $\Sigma_{\io}$ is the polyhedral realization of $B(\ify)$.


Next we determine when $x_{j;i} \equiv 0$ for $j \geq 1$,
$1 \leq i \leq n$.
In the case of $m \geq 1$, $n \leq l \leq 2n-1$
for $\vp_{m;l}(x_{m;1})$ we have 
$$ x_{m+l-n+1;2n-l-1} \geq x_{m+l-n+1;2n-l}$$
and consider the cases that
$(m;l)=(j;2n-1),(j+1;2n-2),\cd,(j+n-1;n)$. Then we have
$$0 \geq x_{j+n;1}\geq x_{j+n;2}\geq \cd \geq x_{j+n;n-1}\geq
x_{j+n;n} \geq 0$$
and $x_{j+n;i} \equiv 0$ for any $j \geq 1$, $1 \leq i \leq n$.
This shows that $\Sigma_{\io}$ is the form of Theorem \ref{c3}.

\end{proof}
\end{thm1}


Next, we give the polyhedral realization $\Sigma_{\io}[\lm]$ for
$\lm := \sum_{k=1}^{n} \lm_k \Lambda_k$, where $\lm_k \in \ZZ_{\geq 0}$, 
$\Lambda_k$ are the fundamental weights. Here we set
\begin{eqnarray*}
\Xi_{\io}^{(1,n-1)} &:=& 
\{ S_{j_k}\cd S_{j_2}  S_{j_1}\xi^{(i)}
(\vec x) : k \geq 0, \; 1 \leq i \leq n-1,\; j_1,\cd,j_k \geq 1 \},\\
\Xi_{\io}^{(n)} &:=& 
\{ S_{j_k}\cd S_{j_2} S_{j_1}\xi^{(n)}
(\vec x) : k \geq 0,\; j_1,\cd,j_k \geq 1 \},\\
\Xi_{\io}^{(1,n-1)}[\lm] &:=& 
\{ \hat S_{j_k}\cd \hat S_{j_2} \hat S_{j_1}\lm^{(i)}
(\vec x) : k \geq 0, \; 1 \leq i \leq n-1,\; j_1,\cd,j_k \geq 1 \},\\
\Xi_{\io}^{(n)}[\lm] &:=& 
\{ \hat S_{j_k}\cd \hat S_{j_2} \hat S_{j_1}\lm^{(n)}
(\vec x) : k \geq 0,\; j_1,\cd,j_k \geq 1 \},\\
\Xi_{\io}[\lm] &:=& \Xi_{\io} \cup \Xi_{\io}^{(1,n-1)}[\lm]\cup
\Xi_{\io}^{(n)}[\lm] \\
	       &=& \{\; \vp_{j;k}(x_{j;1}) :
      j \geq 1, \; 0 \leq k \leq 2n-1 \} \cup 
\{ \hat S_{j_k}\cd \hat S_{j_2} \hat S_{j_1}\lm^{(i)}
(\vec x) : k \geq 0, \; 1 \leq i \leq n-1, \; j_1,\cd,j_k \geq 1 \} \\
&& \; \cup \{ \hat S_{j_k}\cd \hat S_{j_2} \hat S_{j_1}\lm^{(n)}
(\vec x) : k \geq 0,\; j_1,\cd,j_k \geq 1 \}, \\
\Sigma_{\io}[\lm] &:=& \{ \vec x \in \ZZ_{\io}[\lm] : \vp(\vec x)
\geq 0\; \text{ for any }\; \vp \in \Xi_{\io}[\lm] \}.
\end{eqnarray*}
In order to show that $\Sigma_{\io}[\lm]$ is the polyhedral
realization of $B(\lm)$, we give the
explicit forms of $\Xi_{\io}^{(1,n-1)}[\lm]$ and $\Xi_{\io}^{(n)}[\lm]$.  
Since the Dynkin diagrams for $A_n$ and $B_n$ are the same
for $1 \geq i \geq n-1$, 
we give the explicit form of $\Xi_{\io}^{(1,n-1)}$ by the result in
\cite{NZ} as follwos:
\begin{eqnarray}
\Xi_{\io}^{(1,n-1)} = \{ x_{j;i-j} - x_{j;i-j+1} 
: 1 \leq i \leq n-1, \; 1 \leq j \leq i \}.
\label{xiib}
\end{eqnarray}
For giving the explicit form of $\Xi_{\io}^{(n)}$,
we define the ``admissible patterns'' for the integer sequence
$\mu_1,\mu_2,\cd,\mu_k,\cd,\mu_n$ for $1 \leq k \leq n$ as follows:
\begin{eqnarray}
\begin{cases}
 1 \leq \mu_1 \leq n, \\
 0 \leq \mu_2 \leq \mu_1 -1, \\
\q\q\q  \cd, \\
 0 \leq \mu_k \leq \mu_{k-1}-1, \\
\q\q\q  \cd,
\end{cases}
\end{eqnarray}
where if $\mu_k$ does not exist, we define $\mu_k = 0$, and
\begin{eqnarray*}
S^{(\mu_1)}&:=&
\begin{cases}
 id & (\mu_1 = 1), \\
 S_{\mu_1;n-\mu_1} \cd S_{2;n-2}S_{1;n-1} & (\mu_1 \geq2),
\end{cases} \\
S^{(\mu_k)}&:=&
\begin{cases}
  id & (\mu_k=0)\\
  S_{\mu_k+k-2;n-\mu_k+1} \cd S_{k;n-1}S_{k-1;n} & ( \mu_k \geq 1 )
\end{cases}\; \text{ for }\; k \geq 2,\\
\vp^{(\mu)}&:=& S^{(\mu_k)}\cd S^{(\mu_2)}S^{(\mu_1)}\;
\text{ for }\; \mu = (\mu_1,\mu_2,\cd,\mu_k,0,0,\cd).
\end{eqnarray*}
We prepare the symbol $X$ as follows:
\begin{eqnarray}
X := 2x_{1;n-1} - x_{1;n}.
\label{defSX}
\end{eqnarray}
For convenience, we define
\begin{eqnarray}
X_{j;i} :=
\begin{cases}
	2x_{j;i} & \text{ if }\; i \ne n, \\
	x_{j;n}  & \text{ if }\; i = n.
\end{cases} 
\end{eqnarray}

\begin{thm1} Let $\mu = (\mu_1,\mu_2,\cd,\mu_k,0,0,\cd)$ be the
admissible pattern.

\noindent
$(i)$ The forms $\vp^{(\mu)}X$ are given by
\begin{eqnarray}
\vp^{(\mu)}X = \sum_{k=1}^{l}
(X_{\mu_k+k-1;n-\mu_k} - X_{\mu_k+k-1;n-\mu_k+1}),\;
\label{sum_b}
\end{eqnarray}
where
\begin{eqnarray*}
L &:=& \text{max}\{ k : \mu_k \ne 0 \}, \\
l &:=&
\begin{cases}
L &\text{if }\; \mu_L = 1,\\
L +1 &\text{if }\; \mu_L \geq 2.
\end{cases}
\end{eqnarray*}

\noindent
$($ii$)$
$\Xi_{\io}^{(n)}$ 
is the set of all linear forms which are consisted by
$\vp^{(\mu)}X $, where $\mu$ are the admissible patterns.

\begin{proof}
(i) First, we give a remark.
When $\mu_k - \mu_{k+1} =1$, the terms 
$X_{\mu_k+k-1;n-\mu_k+1}$ and $X_{\mu_{k+1}+k;n-\mu_{k+1}}$
in the sum (\ref{sum_b}) are
canceled as $-X_{\mu_k+k-1;n-\mu_k+1} + X_{\mu_{k+1}+k;n-\mu_{k+1}}=0$.
We show the theorem by the induction on 
$|\mu| = \mu_1 + \mu_2 + \cd + \mu_i$
for $\mu = (\mu_1,\mu_2,\cd,\mu_i,0,0,\cd)$.

If $|\mu| = 1$, then $l = 1$ and
the sum of right hand side of (\ref{sum_b}) is
$X_{1;n-1}-X_{1;n} = 2x_{1;n-1}-x_{1;n}$ and equal to $\vp^{(\mu)}X$
by (\ref{defSX}).
We assume that $|\mu| = \mu_1 + \mu_2 + \cd + \mu_i = k-1$ for
$\mu = (\mu_1,\mu_2,\cd,\mu_i,0,0,\cd)$.
We consider the two cases: I) $\mu_i \rightarrow \mu_i +1$,
II) $``\mu_{i+1}=0'' \rightarrow ``\mu_{i+1}=1''$.
\vskip5pt

\noindent
I) $\mu_i \rightarrow \mu_i +1$.

\noindent
In this case, $l = i+1$ and
we note that $1 \leq \mu_i \leq \mu_{i-1} -2$ by admissible
pattern of $\mu$ and the term $X_{\mu_i+i-1;n-\mu_i}$ in the sum
(\ref{sum_b}) is not canceled. We set 
$\mu' := (\mu_1,\mu_2,\cd,\mu_i+1,0,0,\cd)$.
We have
\begin{eqnarray*}
\vp^{(\mu')}X = S_{\mu_i+i-1;n-\mu_i}\vp^{(\mu)}X.
\end{eqnarray*}
The right hand side of (\ref{sum_b}) is
\begin{eqnarray*}
&&S_{\mu_i+i-1;n-\mu_i}(\sum_{k=1}^{l}
  (X_{\mu_k+k-1;n-\mu_k} - X_{\mu_k+k-1;n-\mu_k+1})) \\
=&& \sum_{k=1}^{l}
  (X_{\mu_k+k-1;n-\mu_k} - X_{\mu_k+k-1;n-\mu_k+1}) \\
 &&   - (X_{\mu_i+i-1;n-\mu_i} - X_{\mu_i+i-1;n-\mu_i+1}
     -X_{\mu_i+i;n-\mu_i-1} + X_{\mu_i+i;n-\mu_i}) \\
=&& \vp^{(\mu')}X.
\end{eqnarray*}

\noindent
II) $``\mu_{i+1}=0'' \rightarrow ``\mu_{i+1}=1''$.

\noindent
We set 
$\mu' := (\mu_1,\mu_2,\cd,\mu_i,1,0,\cd)$ $(i.e. \;\mu_{i+1}=1)$
and then $l = i+1$.
By the admissible pattern of $\mu$, we have $0 \leq \mu_{i+1}
\leq \mu_i -1$. This shows that $\mu_i \geq 2$ and
the term $X_{i;n}$ in the sum
(\ref{sum_b}) is not canceled.
We have
\begin{eqnarray*}
\vp^{(\mu')}X = S_{i;n}\vp^{(\mu)}X.
\end{eqnarray*}
The right hand side of (\ref{sum_b}) is
\begin{eqnarray*}
&&S_{i;n}(\sum_{k=1}^{l}
  (X_{\mu_k+k-1;n-\mu_k} - X_{\mu_k+k-1;n-\mu_k+1})) \\
=&& \sum_{k=1}^{l}
  (X_{\mu_k+k-1;n-\mu_k} - X_{\mu_k+k-1;n-\mu_k+1})
   - (X_{i;n} - X_{i+1;n-1} + X_{i+1;n}) \\
=&& \vp^{(\mu')}X.
\end{eqnarray*}

\noindent
(ii) We sufficiently need to show that $\Xi_{\io}^{(n)}$ 
is closed under the
actions of all $S_{j;i}$.
We set for $1 \leq k \leq l+1$,

$\mu =(\mu_1,\mu_2,\cd,\mu_l,0,0,\cd)$, 
$\mu^{+} =(\mu_1,\mu_2,\cd,\mu_k+1,\cd,\mu_l,0,0,\cd)$,

$\mu^{-} =(\mu_1,\mu_2,\cd,\mu_k-1,\cd,\mu_l,0,0,\cd)$.

\noindent
Then we have
\begin{eqnarray*}
S_{j;i}\vp^{(\mu)}X  =
\begin{cases}
\vp^{(\mu^{+})}X 
& \text{ if }\; (j;i)=(\mu_k+k-1;n-\mu_k) \text{ and }\;
\mu_{k-1}-\mu_k \ne 1, \\
\vp^{(\mu^{-})}X 
& \text{ if }\; (j;i)=(\mu_k+k-1;n-\mu_k+1) \text{ and }\;
\mu_{k}-\mu_{k+1} \ne 1, \\
\vp^{(\mu)}X & \text{ other wise. }
\end{cases}
\end{eqnarray*}

\end{proof}
\end{thm1}
Therefore, $\io$ satisfies the strict positivity assumption
by the explicit forms of $\Xi_{\io}^{(1,n-1)}$
and $\Xi_{\io}^{(n)}$. 
By the Lemma \ref{lmxi}, this shows that
\begin{eqnarray*}
\Xi_{\io}^{(1,n-1)}[\lm] &=&
 \{ \lm_i + x_{j;i-j} - x_{j;i-j+1} 
 : 1 \leq i \leq n-1, \; 1 \leq j \leq i \}, \\
\Xi_{\io}^{(n)}[\lm] &=&
 \{ \lm_n + \sum_{k=1}^{l}
 (X_{\mu_k+k-1;n-\mu_k} - X_{\mu_k+k-1;n-\mu_k+1}) :
 \mu \; \text{ are the admissible patterns }\}
\end{eqnarray*}
by (\ref{xiib}), (\ref{sum_b}) and
$\Sigma_{\io}[\lm]$ is the polyhedral realization of $B(\lm)$.


\subsection{$C_n$ case}

We consider the case of type $C_n$. First, we give the
polyhedral realization $\Sigma_{\io}$ of $B(\ify)$.
We define for $j \geq 1,\, 0 \leq k \leq
2n-1$
\begin{eqnarray*}
\vp_{j;k}:=
\begin{cases}
  id & (k=0),\\
  S_{j;k}S_{j;k-1} \cd S_{j;2}S_{j;1} & (1 \leq k \leq n-1),\\
  S_{j+k-n;2n-k}S_{j+k-1-n;2n-k+1} \cd S_{j+1;n-1}S_{j;n}\vp_{j;n-1}
 & (n \leq k \leq 2n-1).
\end{cases}
\label{b_phi_s}
\end{eqnarray*}

\begin{lem1}
\begin{eqnarray*}
\vp_{j;k}(x_{j;1})=
\begin{cases}
x_{j;k+1}-x_{j+1;k} & (0 \leq k \leq n-2),\\
2x_{j;n}-x_{j+1;n-1} & (k = n-1), \\
x_{j+1;n-1}-2x_{j+1;n} & (k = n), \\
x_{j+k-n+1;2n-k-1}-x_{j+k-n+1;2n-k} & (n+1 \leq k \leq 2n-1).
\end{cases}
\label{b_phi_x}
\end{eqnarray*}

\begin{proof}
By the induction on $k$. The case of $k=0$ is trivial.\\

\noindent
Case I) $1 \leq k \leq n-2$.\\

\noindent
If $k=1$, we have
\begin{eqnarray*}
S_{j;1}(x_{j;1}) &=& x_{j;1} - (x_{j;1}-x_{j;2}+x_{j+1;1}) \\
		 &=& x_{j;2}-x_{j+1;1} \\
		 &=& \vp_{j;1}(x_{j;1}).
\end{eqnarray*}

\noindent
If $1 < k \leq n-2$, we assume that 
$\vp_{j;k-1}(x_{j;1}) = x_{j;k}-x_{j+1;k-1}$. Then we have
\begin{eqnarray*}
S_{j;k}\vp_{j;k-1}(x_{j;1}) &=& S_{j;k}(x_{j;k}-x_{j+1;k-1}) \\ 
        &=& x_{j;k}-x_{j+1;k-1}-(x_{j;k}-x_{j;k+1}-x_{j+1;k-1}+x_{j+1;k}) \\
	&=& x_{j;k+1}-x_{j+1;k} \\
	&=& \vp_{j;k}(x_{j;1}).
\end{eqnarray*}

\noindent
Case II) $k = n-1,n$.\\

\noindent
If $k=n-1$, using the result of Case I), we have
\begin{eqnarray*}
S_{j;n-1}\vp_{j;n-2}(x_{j;1}) &=& S_{j;n-1}(x_{j;n-1}-x_{j+1;n-2}) \\
  &=& x_{j;n-1}-x_{j+1;n-2}-(x_{j;n-1}-2x_{j;n}-x_{j+1;n-2}+x_{j+1;n-1}) \\
	&=& 2x_{j;n}-x_{j+1;n-1} \\
	&=& \vp_{j;n-1}(x_{j;1}),
\end{eqnarray*}
and then if $k=n$, 
\begin{eqnarray*}
S_{j;n}\vp_{j;n-1}(x_{j;1}) &=& S_{j;n}(2x_{j;n}-x_{j+1;n-1}) \\
  &=& 2x_{j;n}-x_{j+1;n-1}-2(x_{j;n}-x_{j+1;n-1}+x_{j+1;n}) \\
	&=& x_{j+1;n-1}-2x_{j+1;n} \\
	&=& \vp_{j;n}(x_{j;1}).
\end{eqnarray*}

\noindent
Case III) $n+1 \leq k \leq 2n-1$.\\

\noindent
If $k=n+1$, using the result of the case of $k=n$, we have
\begin{eqnarray*}
S_{j+1;n-1}\vp_{j;n}(x_{j;1}) &=& S_{j+1;n-1}(x_{j+1;n-1}-2x_{j+1;n}) \\
&=& x_{j+1;n-1}-2x_{j+1;n}-(x_{j+1;n-1}-2x_{j+1;n}-x_{j+2;n-2}+x_{j+2;n-1}) \\
	&=& x_{j+2;n-2}-x_{j+2;n-1} \\
	&=& \vp_{j;n+1}(x_{j;1}).
\end{eqnarray*}

\noindent
If $n+1 < k \leq 2n-1$, we assume that 
$\vp_{j;k-1}(x_{j;1}) = x_{j+k-n;2n-k}-x_{j+k-n;2n-k+1}$. Then we have 
\begin{eqnarray*}
S_{j+k-n;2n-k}\vp_{j;k-1}(x_{j;1}) &=& S_{j+k-n;2n-k}
(x_{j+k-n;2n-k}-x_{j+k-n;2n-k+1}) \\ 
&=& x_{j+k-n;2n-k}-x_{j+k-n;2n-k+1} \\
&&-(x_{j+k-n;2n-k}-x_{j+k-n;2n-k+1}-x_{j+k-n+1;2n-k-1}+x_{j+k-n+1;2n-k}) \\
&=& x_{j+k-n+1;2n-k-1}-x_{j+k-n+1;2n-k} \\
&=& \vp_{j;k}(x_{j;1}).
\end{eqnarray*}
\end{proof}
\label{b1}
\end{lem1}

\begin{lem1}
$\vp_{j;k}(x_{j;1})$ is closed
under the actions of all transformations $S_{m;l}$ for any
$m \geq 1, l \geq 1$.
\begin{proof}
Using the definition of $\vp_{j;k}$ (\ref{b_phi_s}) and the formula
(\ref{b_phi_x}) in Lemma \ref{b1},
if $k=0$, $\vp_{j;0} = id$ and we have
\begin{eqnarray*}
S_{m;l}(\vp_{j;0}(x_{j;1}))=
\begin{cases}
\vp_{j;1}(x_{j;1}) & if \;\;  (m;l)=(j;1),\\
x_{j;1} & other \; wise.
\end{cases}
\end{eqnarray*}

If $1 \leq k \leq n-1$, we have
\begin{eqnarray*}
S_{m;l}(\vp_{j;k}(x_{j;1}))=
\begin{cases}
\vp_{j;k-1}(x_{j;1}) & if \;\;  (m;l)=(j;k),\\
\vp_{j;k+1}(x_{j;1}) & if \;\;  (m;l)=(j;k+1),\\
\vp_{j;k}(x_{j;1}) & other \; wise.
\end{cases}
\end{eqnarray*}

If $n \leq k \leq 2n-1$, we have
\begin{eqnarray*}
S_{m;l}(\vp_{j;k}(x_{j;1}))=
\begin{cases}
\vp_{j;k-1}(x_{j;1}) & if \;\;  (m;l)=(j+k-n;2n-k),\\
\vp_{j;k+1}(x_{j;1}) & if \;\;  (m;l)=(j+k-n+1;2n-k-1),\\
\vp_{j;k}(x_{j;1}) & other \; wise.
\end{cases}
\end{eqnarray*}

\end{proof}
\end{lem1}

Now, we define
\begin{eqnarray*}
\Xi_{\io} &:=& \{\; \vp_{j;k}(x_{j;1}) :
j \geq 1,\; 0 \leq k \leq 2n-1 \}, \\
\Sigma_{\io} &:=& \{\; \vec x \in \ZZ^{\ify}_{\io} \;|\;
\vp(\vec x) \geq 0 \text{ for any } \vp \in \Xi_{\io} \;\}.
\end{eqnarray*}

\begin{thm1}
Polyhedral realization $\Sigma_{\io}$
of $B(\ify)$ for type $C_n$ is described as
below:
\begin{eqnarray*}
&&x_{j;i} = 0 \q \text{ for }\; j, i \notin [1,n], \\
&&x_{1;i}\geq x_{2;i-1}\geq \cd \geq x_{i;1} \geq 0
 \q \text{ for }\; 1 \leq i \leq n-1,\\
&&2x_{j;n}\geq x_{j+1;n-1}\geq \cd \geq x_{n;j}\geq 0
\q \text{ for }\; 1 \leq j \leq n-1,\\
&&x_{j;n-j+1} \geq x_{j;n-j+2} \geq \cd \geq x_{j;n-1} \geq
2x_{j;n} \geq 0
\q \text{ for }\; 2 \leq j \leq n.
\end{eqnarray*}
\begin{proof}

We shall show that $\Sigma_{\io}$ is the polyhedral realization
of $B(\ify)$ and check the conditions of Theorem \ref{nz_modi}.
$\Xi_{\io}$ is closed under the action of $S_k$ by Lemma \ref{b2}.
The coefficients of $x_{1;i}$ $(i = 1,2,\cd,n)$ are
positive for $\vp \in \Xi_{\io}$ by Lemma $\ref{b1}$.
This shows that $\io$ satisfies the positivity assumption.

We shall show that all entries of $\vec x \in \Sigma_{\io}$
are nonnegative. In the case of $m \geq 1$, $0 \leq l \leq n-2$ for
$\vp_{m;l}(x_{m;1})$, we have 
$$x_{m;l+1} \geq x_{m+1;l}$$
and we consider the cases that 
$(m;l)=(j;i-1),(j+1;i-2),\cd,(j+i-2;1),(j+i-1;0)$, then we obtain
\begin{eqnarray}
x_{j;i}\geq x_{j+1;i-1} \geq \cd \geq 
x_{j+i-2;2}\geq x_{j+i-1;1} \geq 0.
\label{ccc}
\end{eqnarray}
This shows that $x_{j;i} \geq 0$ for any $j \geq 1$, $1 \leq i \leq
n-1$.
Similarly, for any $j \geq 1$ and $l=n-1$, we have 
$$2x_{j;n} \geq x_{j+1;n-1} \;(\geq 0).$$
This shows $x_{j;n} \geq 0$ for any $j \geq 1$.
Therefore, $\Sigma_{\io}$ is the polyhedral realization of $B(\ify)$.

On the other hand, 
in the case of $(m;l)=(j+n-1;n)$ for $\vp_{m;l}(x_{m;1})$, we obtain
$$ x_{j+n;n-1} \geq 2x_{j+n;n}$$
and in the case of $m \geq 1$, $n+1 \leq l \leq 2n-1,$
we have 
$$ x_{m+l-n+1;2n-l-1} \geq x_{m+l-n+1;2n-l}.$$
We consider the cases that
$(m;l)=(j;2n-1),(j+1;2n-2),\cd,(j+n-2;n+1)$, then we have
$$0 \geq x_{j+n;1}\geq x_{j+n;2}\geq \cd \geq 
x_{j+n;n-2} \geq x_{j+n;n-1}\geq 2x_{j+n;n} \geq 0.$$
This shows that $x_{j+n;i} \equiv 0$ for any $j \geq 1$, $1 \leq i \leq n$.

\end{proof}
\end{thm1}


Next, we give the polyhedral realization $\Sigma_{\io}[\lm]$
of $B(\lm)$.
We prepare the following symbols:
\begin{eqnarray}
X &:=& x_{1;n-1} - x_{1;n}, \\
X_{j;i} &:=&
\begin{cases}
	x_{j;i} & \text{ if }\; i \ne n, \\
	2x_{j;n}  & \text{ if }\; i = n.
\end{cases} 
\end{eqnarray}
Then, we can show that the polyhedral realization $\Sigma_{\io}[\lm]$
is given by the replacing $x_{j;i}$ of $B_n$ with $X_{j;i}$.


\subsection{$D_n$ case}

We consider the case of type $D_n$. First, we give the
polyhedral realization $\Sigma_{\io}$ of $B(\ify)$.
We define for $j \geq 1,\, 0 \leq k \leq
2n-2$  
\begin{align*}
\vp_{j;k}:= &
\begin{cases}
  id & (k = 0),\\
  S_{j;k}S_{j;k-1} \cd S_{j;2}S_{j;1} & (0 \leq k \leq n),\\
  S_{j+k-n;2n-k-1}S_{j+k-1-n;2n-k} \cd S_{j+2;n-3}S_{j+1;n-2}\vp_{j;n}
 & (n+1 \leq k \leq 2n-2).
\end{cases}
\\
\vp'_{j;k}:= &
\begin{cases}
  id & (k = 0),\\
  S_{j;k}S_{j;k-1} \cd S_{j;2}S_{j;1} & (1 \leq k \leq n-2),\\
  S_{j;n}\vp'_{j;n-2} & (k = n-1),\\
  S_{j;n-2}\vp'_{j;n-1} & (k = n),\\
  S_{j+k-n;2n-k-1}S_{j+k-1-n;2n-k} \cd S_{j+2;n-3}S_{j+1;n-2}\vp'_{j;n}
 & (n+1 \leq k \leq 2n-2).
\end{cases}
\end{align*}

\begin{lem1}
\begin{align*}
\vp_{j;k}(x_{j;1})= &
\begin{cases}
x_{j;k+1}-x_{j+1;k} & (0 \leq k \leq n-3),\\
x_{j;n-1}+x_{j;n}-x_{j+1;n-2} & (k = n-2), \\
x_{j;n}-x_{j+1;n-1} & (k = n-1), \\
x_{j+1;n-2}-x_{j+1;n-1}-x_{j+1;n} & (k = n), \\
x_{j+k-n+1;2n-k-2}-x_{j+k-n+1;2n-k-1} & (n+1 \leq k \leq 2n-2).
\end{cases}
\\
\vp'_{j;k}(x_{j;1})= &
\begin{cases}
x_{j;n-1}-x_{j+1;n} & (k = n-1), \\
\vp_{j;k}(x_{j;1}) & other \; wise.
\end{cases}
\end{align*}
\begin{proof}
We shall prove the case $\vp_{j;k}$ by the induction on $k$
since the other case
can be proved by the same argument.
The case of $k=0$ is trivial.\\

\noindent
Case I) $1 \leq k \leq n-3$.\\

\noindent
If $k=1$, we have
\begin{eqnarray*}
S_{j;1}(x_{j;1}) &=& x_{j;1} - (x_{j;1}-x_{j;2}+x_{j+1;1}) \\
		 &=& x_{j;2}-x_{j+1;1} \\
		 &=& \vp_{j;1}(x_{j;1}).
\end{eqnarray*}

\noindent
If $1 < k \leq n-3$, we assume that 
$\vp_{j;k-1}(x_{j;1}) = x_{j;k}-x_{j+1;k-1}$. Then we have
\begin{eqnarray*}
S_{j;k}\vp_{j;k-1}(x_{j;1}) &=& S_{j;k}(x_{j;k}-x_{j+1;k-1}) \\ 
        &=& x_{j;k}-x_{j+1;k-1}-(x_{j;k}-x_{j;k+1}-x_{j+1;k-1}+x_{j+1;k}) \\
	&=& x_{j;k+1}-x_{j+1;k} \\
	&=& \vp_{j;k}(x_{j;1}).
\end{eqnarray*}

\noindent
Case II) $k = n-2,n-1,n$.\\

\noindent
If $k=n-2$, using the result of Case I), we have
\begin{eqnarray*}
S_{j;n-2}\vp_{j;n-3}(x_{j;1}) &=& S_{j;n-2}(x_{j;n-2}-x_{j+1;n-3}) \\
  &=& x_{j;n-2}-x_{j+1;n-3}-(x_{j;n-2}-x_{j;n-1}-x_{j;n}
-x_{j+1;n-3}+x_{j+1;n-2}) \\
	&=& x_{j;n-1}+x_{j;n}-x_{j+1;n-2} \\
	&=& \vp_{j;n-2}(x_{j;1}),
\end{eqnarray*}
and then if $k=n-1$, 
\begin{eqnarray*}
S_{j;n-1}\vp_{j;n-2}(x_{j;1}) &=& S_{j;n-1}(x_{j;n-1}+x_{j;n}-x_{j+1;n-2}) \\
  &=& x_{j;n-1}+x_{j;n}-x_{j+1;n-2}
-(x_{j;n-1}-x_{j+1;n-2}+x_{j+1;n-1}) \\
	&=& x_{j;n}-x_{j+1;n-1} \\
	&=& \vp_{j;n-1}(x_{j;1}).
\end{eqnarray*}

If $k=n$, then we have

\begin{eqnarray*}
S_{j;n}\vp_{j;n-1}(x_{j;1}) &=& S_{j;n}(x_{j;n}-x_{j+1;n-1}) \\
  &=& x_{j;n}-x_{j+1;n-1}
-(x_{j;n}-x_{j+1;n-2}+x_{j+1;n}) \\
	&=& x_{j+1;n-2}-x_{j+1;n-1}-x_{j+1;n} \\
	&=& \vp_{j;n}(x_{j;1}).
\end{eqnarray*}

\noindent
Case III) $n+1 \leq k \leq 2n-2$.\\

\noindent
If $k=n+1$, using the result of the case of $k=n$, we have
\begin{eqnarray*}
S_{j+1;n-2}\vp_{j;n}(x_{j;1}) &=& S_{j+1;n-2}
(x_{j+1;n-2}-x_{j+1;n-1}-x_{j+1;n}) \\
&=& x_{j+1;n-2}-x_{j+1;n-1}-x_{j+1;n} \\
&&-(x_{j+1;n-2}-x_{j+1;n-1}-x_{j+1;n}-x_{j+2;n-3}+x_{j+2;n-2}) \\
	&=& x_{j+2;n-3}-x_{j+2;n-2} \\
	&=& \vp_{j;n+1}(x_{j;1}).
\end{eqnarray*}

\noindent
If $n+1 < k \leq 2n-2$, we assume that
$\vp_{j;k-1}(x_{j;1}) = x_{j+k-n;2n-k-1}-x_{j+k-n;2n-k}$.
Then we have

\begin{eqnarray*}
S_{j+k-n;2n-k-1}\vp_{j;k-1}(x_{j;1}) &=& S_{j+k-n;2n-k-1}
(x_{j+k-n;2n-k-1}-x_{j+k-n;2n-k}) \\
&=& x_{j+k-n;2n-k-1}-x_{j+k-n;2n-k} \\
&&-(x_{j+k-n;2n-k-1}-x_{j+k-n;2n-k}-x_{j+k-n+1;2n-k-2}
+x_{j+k-n+1;2n-k-1}) \\
	&=&  x_{j+k-n+1;2n-k-2}-x_{j+k-n+1;2n-k-1})\\
	&=& \vp_{j;k}(x_{j;1}).
\end{eqnarray*}
\end{proof}
\label{d1}
\end{lem1}

\begin{lem1}
The set of $\vp_{j;k}(x_{j;1})$ and $\vp'_{j;k}(x_{j;1})$ is closed
under the actions of all transformations $S_{m;l}$
for any $m \geq 1, l \geq 1$.
\begin{proof}
First, we consider the action of $S_{m;l}$ to $\vp_{j;k}(x_{j;1})$.
If $k=0$, $\vp_{j;0}=id$ and we have
\begin{eqnarray*}
S_{m;l}(\vp_{j;0}(x_{j;1}))=
\begin{cases}
\vp_{j;1}(x_{j;1}) & if \;\;  (m;l)=(j;1),\\
x_{j;1} & other \;wise.
\end{cases}
\end{eqnarray*}

If $1 \leq k \leq n-3$, we have
\begin{eqnarray*}
S_{m;l}(\vp_{j;k}(x_{j;1}))=
\begin{cases}
\vp_{j;k-1}(x_{j;1}) & if \;\;  (m;l)=(j;k),\\
\vp_{j;k+1}(x_{j;1}) & if \;\;  (m;l)=(j;k+1),\\
\vp_{j;k}(x_{j;1}) & other \;wise.
\end{cases}
\end{eqnarray*}

If $k = n-2$, we have
\begin{eqnarray}
S_{m;l}(\vp_{j;n-2}(x_{j;1}))=
\begin{cases}
\vp_{j;n-3}(x_{j;1}) & if \;\;  (m;l)=(j+1;n-2),\\
\vp_{j;n-1}(x_{j;1}) & if \;\;  (m;l)=(j;n-1),\\
\vp'_{j;n-1}(x_{j;1}) & if \;\;  (m;l)=(j;n),\\
\vp_{j;n-2}(x_{j;1}) & other \;wise.
\end{cases}
\label{vp'n-2}
\end{eqnarray}

If $k = n-1$, we have
\begin{eqnarray*}
S_{m;l}(\vp_{j;n-1}(x_{j;1}))=
\begin{cases}
\vp_{j;n-2}(x_{j;1}) & if \;\;  (m;l)=(j+1;n-1),\\
\vp_{j;n}(x_{j;1}) & if \;\;  (m;l)=(j;n),\\
\vp_{j;n-1}(x_{j;1}) & other \;wise.
\end{cases}
\end{eqnarray*}

If $k = n$, we have
\begin{eqnarray}
S_{m;l}(\vp_{j;n}(x_{j;1}))=
\begin{cases}
\vp_{j;n-1}(x_{j;1}) & if \;\;  (m;l)=(j+1;n),\\
\vp'_{j;n-1}(x_{j;1}) & if \;\;  (m;l)=(j+1;n-1),\\
\vp_{j;n+1}(x_{j;1}) & if \;\;  (m;l)=(j+1;n-2),\\
\vp_{j;n}(x_{j;1}) & other \;wise.
\end{cases}
\label{vp'n}
\end{eqnarray}

If $n+1 \leq k \leq 2n-2$, we have
\begin{eqnarray*}
S_{m;l}(\vp_{j;k}(x_{j;1}))=
\begin{cases}
\vp_{j;k-1}(x_{j;1}) & if \;\;  (m;l)=(j+k-n+1;2n-k-1),\\
\vp_{j;k+1}(x_{j;1}) & if \;\;  (m;l)=(j+k-n+1;2n-k-2),\\
\vp_{j;k}(x_{j;1}) & other \;wise.
\end{cases}
\end{eqnarray*}

Next, we consider the action of $S_{m;l}$ to $\vp'_{j;k}(x_{j;1})$.
Since the difference of $\vp_{j;k}(x_{j;1})$ and
$\vp'_{j;k}(x_{j;1})$ is only the case of $k = n-1$ and
by (\ref{vp'n-2}), (\ref{vp'n}),
we consider the case of $k = n-1$.

\begin{eqnarray*}
S_{m;l}(\vp'_{j;n-1}(x_{j;1}))=
\begin{cases}
\vp'_{j;n-2}(x_{j;1})\;(=\vp_{j;n-2}(x_{j;1})) & if \;\;  (m;l)=(j+1;n),\\
\vp'_{j;n}(x_{j;1})\;(=\vp_{j;n}(x_{j;1})) & if \;\;  (m;l)=(j;n-1),\\
\vp'_{j;n-1}(x_{j;1}) & other \;wise.
\end{cases}
\end{eqnarray*}

\end{proof}
\label{d2}
\end{lem1}

Now, we define
\begin{eqnarray*}
&\Xi'_{\io}& := \;\; \{\; \vp_{j;k}(x_{j;1}) :
j \geq 1, \;0 \leq k \leq 2n-2 \}, \\
&\Xi''_{\io}& := \;\;\{ x_{j;n-1} : j \geq 1\}
\cup \{ x_{j;n} : j \geq 1\}, \\
&\Xi_{\io}& := \;\;\Xi'_{\io}\cup\Xi''_{\io}, \\
&\Sigma_{\io}& := \;\;\{\; \vec x \in \ZZ^{\ify}_{\io} :
\vp(\vec x) \geq 0 \text{ for any } \vp \in \Xi_{\io}\;\}.
\end{eqnarray*}

\begin{thm1}
Polyhedral realization $\Sigma_{\io}$
of $B(\ify)$ for type $D_n$ is described as
below:
\label{d3}
\begin{eqnarray*}
&&x_{j;i} = 0 \q \text{ for }\; j \notin [1,n-1]\;\;or\;\; i \notin [1,n],\\
&&x_{1;i}\geq x_{2;i-1}\geq \cd \geq x_{i;1} \geq 0
 \q \text{ for }\; 1 \leq i \leq n-2,\\
&&x_{j;n-1}+x_{j;n}\geq x_{j+1;n-2}\geq x_{j+2;n-3}
\cd \geq x_{n-1;j}\geq 0
\q \text{ for }\; 1 \leq j \leq n-2,\\
&&x_{j;n-j} \geq x_{j;n-j+1} \geq \cd \geq x_{j;n-2} \geq
x_{j;n-1}+x_{j;n} \geq 0
\q \text{ for }\; 2 \leq j \leq n-1,\\
&&
\begin{cases}
x_{1;n-1} \geq x_{2;n} \geq x_{3;n-1} \geq x_{4;n} \geq \cd
\geq x_{n-1;n} \geq 0 \\
x_{1;n} \geq x_{2;n-1} \geq x_{3;n} \geq x_{4;n-1} \geq \cd
\geq x_{n-1;n-1} \geq 0
\end{cases} ( n: odd)
\\
&&
\begin{pmatrix}
\;
\begin{cases}
x_{1;n-1} \geq x_{2;n} \geq x_{3;n-1} \geq x_{4;n} \geq \cd
\geq x_{n-1;n-1} \geq 0 \\
x_{1;n} \geq x_{2;n-1} \geq x_{3;n} \geq x_{4;n-1} \geq \cd
\geq x_{n-1;n} \geq 0
\end{cases} (resp. \;\; n: even)
\;
\end{pmatrix}.
\end{eqnarray*}
\begin{proof}

Since $\Xi_{\io}$ is not closed under the action of $S_k$, we
shall show that

\noindent
$(i)$ $\Sigma_{\io}$ is closed under the action of $\fit$,

\noindent
$(ii)$ $\Sigma_{\io}$ is closed under the action of $\eit$,

\noindent
$(iii)$ all the entries of $\vec x \in \Sigma_{\io}$ are nonnegative

\noindent
(These shows that $\Sigma_{\io}$ is the polyhedral realization of
$B(\ify)$),

\noindent
$(iv)$ $\Sigma_{\io}$ is the form of Theorem \ref{d3}.
\vskip5pt

\noindent
$(i)$ We show that $\Sigma_{\io}$ is closed under the action of $\fit$.

Note that $\io$ satisfies the positivity assumption by Lemma \ref{d2}.
Let $\vec x = (\cd,x_2,x_1) \in \Sigma_{\io}$ and $i \in I$, and
suppose that $\fit \vec x = (\cd,x_k +1, \cd, x_2,x_1)$ for $i_k = i$.
We need to show that
$$\vp(\fit \vec x) \geq 0$$
for any $\vp = \sum\vp_jx_j \in \Xi_{\io}$.
First, we consider the case of $\vp \in \Xi_{\io}'$.
Since $\vp(\fit\vec x) = \vp(\vec x) + \vp_k$, it is enough to
consider the case when $\vp_k < 0$.
Since $\io$ satisfies the positivity assumption, 
we have $k^{(-)} \geq 1$. By (\ref{action-f}), we have
$\sigma_k(\vec x) > \sigma_{k^{(-)}}(\vec x)$ and by (\ref{beta}),
we conclude that
$$\beta_{k^{(-)}}(\vec x)=\sigma_{k^{(-)}}(\vec x) 
-\sigma_{k}(\vec x) \leq -1.$$

It follows that
\begin{eqnarray}
\vp(\fit \vec x) &=& \vp(\vec x) + \vp_k  \nn\\
		 &\geq& \vp(\vec x)-1 \cdot \beta_{k^{(-)}}(\vec x) \nn\\
		 &=& (S_k \vp)(\vec x) \\
		 &\geq& 0
\label{dddd}
\end{eqnarray}
since $S_k \vp \in \Xi_{\io}'$.
This shows that $\Sigma_{\io}$ is closed under the action of
$\fit$ for $\vp \in \Xi_{\io}'$. 

Next, we consider the case of $\vp = x_{j;n-1} \in \Xi_{\io}''$. We have
\begin{eqnarray*}
\vp(\fit \vec x) &=& x_{j;n-1} +1 \\
		 &\geq& 0.
\end{eqnarray*}
The case of $\vp = x_{j;n} \in \Xi''_{\io}$ can be poved by the same
argument.
\vskip5pt

\noindent
$(ii)$ We show that $\Sigma_{\io}$ is closed under the action of $\eit$.

We need to show that
$$\vp(\eit \vec x) \geq 0$$
for any $\vp \in \Xi_{\io}$. First, we consider the case of
$\vp \in \Xi_{\io}'$.

Since $\vp(\eit\vec x) = \vp(\vec x) - \vp_k$, it is enough to
consider the case when $\vp_k > 0$.
By (\ref{action-e}),
$\sigma_k(\vec x) > \sigma_{k^{(+)}}(\vec x)$ and by (\ref{beta}),
we conclude that
$$\beta_{k}(\vec x)=\sigma_{k}(\vec x) 
-\sigma_{k^{(+)}}(\vec x) \geq 1.$$

It follows that
\begin{eqnarray*}
\vp(\eit \vec x) &=& \vp(\vec x) - \vp_k \\
		 &\geq& \vp(\vec x)-1 \cdot \beta_{k}(\vec x) \\
		 &=& (S_k \vp)(\vec x) \\
		 &\geq& 0
\end{eqnarray*}
since $S_k \vp \in \Xi_{\io}'$.
This shows that $\Sigma_{\io}$ is closed under the action of
$\eit$ for $\vp \in \Xi_{\io}'$.

Next, we consider the case of $\vp = x_{j;n-1} \in \Xi_{\io}''$.
It follows that
\begin{eqnarray*}
\vp(\eit \vec x) &=& x_{j;n-1}(\vec x) - 1 \\
		 &\geq& x_{j;n-1}(\vec x)-1 \cdot \beta_{j;n-1}(\vec x) \\
		 &=& (S_{j;n-1}x_{j;n-1})(\vec x) \\
		 &=& (x_{j+1;n-2} - x_{j+1;n-1})(\vec x) \\
		 &\geq& x_{j+1;n}(\vec x) \q (by \; Lemma\; \ref{d1}) \\
		 &\geq& 0 \q (by \; \;x_{j+1;n} \in \Xi''_{\io}).
\end{eqnarray*}
The case of $\vp = x_{j;n} \in \Xi''_{\io}$ can be poved by the same
argument.
\vskip5pt

\noindent
$(iii)$ We show that all the entries of $\vec x \in \Sigma_{\io}$ are
nonnegative. 

In the case of $m \geq 1$, $0 \leq l \leq n-3$ for
$\vp_{m;l}(x_{m;1})$, we have 
$$x_{m;l+1} \geq x_{m+1;l}, $$
and we consider the cases that 
$(m;l)=(j;i-1),(j+1;i-2),\cd,(j+i-2;1),(j+i-1;0)$, then we obtain
\begin{eqnarray}
x_{j;i}\geq x_{j+1;i-1} \geq \cd \geq 
x_{j+i-2;2}\geq x_{j+i-1;1} \geq 0.
\label{b2}
\end{eqnarray}
This shows that $x_{j;i} \geq 0$ for any $j \geq 1$, $1 \leq i \leq
n-2$.
By the definition of $\Xi''_{\io}$, we have $x_{j;n-1}$ and
$x_{j;n} \geq 0$ for $j \geq 1$
and these shows $x_{j;n} \geq 0$ for any $j \geq 1$.
\vskip5pt

\noindent
$(iv)$ We determine when $x_{j;i} \equiv 0$ for $j\geq1$, 
$1 \leq i\leq n$.

In the case of $m \geq 1$, $n+1 \leq l \leq 2n-2$
for $\vp_{m;l}(x_{m;1})$, we have 
$$ x_{m+l-n+1;2n-l-2} \geq x_{m+l-n+1;2n-l-1}$$
and in the case of $l = n$,
$$x_{m+1;n-2} \geq x_{m+1;n-1}+x_{m+1;n} \geq 0 \;\;(\because
\; x_{m+1;n-1} \geq 0 \text{ and }x_{m+1;n}\geq 0).$$
We consider the cases that
$(m;l)=(j;2n-2),(j+1;2n-3),\cd,(j+n-2;n)$, then we have
$$0 \geq x_{j+n-1;1}\geq x_{j+n-1;2}\geq \cd \geq 
x_{j+n-1;n-3} \geq x_{j+n-1;n-2}\geq x_{j+n-1;n-1}+x_{j+n-1;n} \geq 0.$$
Combining $x_{j+n-1;n-1} \geq 0, x_{j+n-1;n}\geq 0$
this shows that $x_{j+n-1;i} \equiv 0$ for any $j \geq 1$, $1 \leq i \leq n$.


\end{proof}
\end{thm1}


Next, we give the polyhedral realization $\Sigma_{\io}[\lm]$
of $B(\lm)$ for
$\lm := \sum_{k=1}^{n} \lm_k \Lambda_k$, where $\lm_k \in \ZZ_{\geq 0}$, 
$\Lambda_k$ are the fundamental weights. Here we set
\begin{eqnarray*}
\Xi_{\io}^{(1,n-2)} &:=& 
\{ S_{j_k}\cd S_{j_2}  S_{j_1}\xi^{(i)}
(\vec x) : k \geq 0, \; 1 \leq i \leq n-2,\; j_1,\cd,j_k \geq 1 \},\\
\Xi_{\io}^{(n-1)} &:=& 
\{ S_{j_k}\cd S_{j_2} S_{j_1}\xi^{(n-1)}
(\vec x) : k \geq 0,\; j_1,\cd,j_k \geq 1 \},\\
\Xi_{\io}^{(n)} &:=& 
\{ S_{j_k}\cd S_{j_2} S_{j_1}\xi^{(n)}
(\vec x) : k \geq 0,\; j_1,\cd,j_k \geq 1 \},\\
\Xi_{\io}^{(1,n-2)}[\lm] &:=& 
\{ \hat S_{j_k}\cd \hat S_{j_2} \hat S_{j_1}\lm^{(i)}
(\vec x) : k \geq 0, \; 1 \leq i \leq n-2,\; j_1,\cd,j_k \geq 1 \},\\
\Xi_{\io}^{(n-1)}[\lm] &:=& 
\{ \hat S_{j_k}\cd \hat S_{j_2} \hat S_{j_1}\lm^{(n-1)}
(\vec x) : k \geq 0,\; j_1,\cd,j_k \geq 1 \},\\
\Xi_{\io}^{(n)}[\lm] &:=& 
\{ \hat S_{j_k}\cd \hat S_{j_2} \hat S_{j_1}\lm^{(n)}
(\vec x) : k \geq 0,\; j_1,\cd,j_k \geq 1 \},\\
\Xi_{\io}[\lm] &:=& \Xi_{\io} \cup \Xi_{\io}^{(1,n-2)}[\lm]\cup
\Xi_{\io}^{(n-1)}[\lm]\cup\Xi_{\io}^{(n)}[\lm] \\
	       &=& \{\; \vp_{j;k}(x_{j;1}) :
      j \geq 1, \; 0 \leq k \leq 2n-1 \} \cup 
\{ \hat S_{j_k}\cd \hat S_{j_2} \hat S_{j_1}\lm^{(i)}
(\vec x) : k \geq 0, \; 1 \leq i \leq n-2, \; j_1,\cd,j_k \geq 1 \} \\
&& \; \cup \{ \hat S_{j_k}\cd \hat S_{j_2} \hat S_{j_1}\lm^{(n-1)}
(\vec x) : k \geq 0,\; j_1,\cd,j_k \geq 1 \}\cup
\{ \hat S_{j_k}\cd \hat S_{j_2} \hat S_{j_1}\lm^{(n)}
(\vec x) : k \geq 0,\; j_1,\cd,j_k \geq 1 \}, \\
\Sigma_{\io}[\lm] &:=& \{ \vec x \in \ZZ_{\io}[\lm] : \vp(\vec x)
\geq 0\; \text{ for any }\; \vp \in \Xi_{\io}[\lm] \}.
\end{eqnarray*}
In order to show that $\Sigma_{\io}[\lm]$ is the polyhedral
realization of $B(\lm)$, we give the
explicit forms of $\Xi_{\io}^{(1,n-2)}[\lm]$, $\Xi_{\io}^{(n-1)}[\lm]$ and
$\Xi_{\io}^{(n)}[\lm]$.  
Since the Dynkin diagrams for $A_n$ and $D_n$ are the same
for $1 \leq i \leq n-2$,
we give the explicit form of $\Xi_{\io}^{(1,n-2)}$ 
 as follwos:
\begin{eqnarray}
\Xi_{\io}^{(1,n-2)} = \{ x_{j;i-j} - x_{j;i-j+1} 
: 1 \leq i \leq n-2, \; 1 \leq j \leq i \}.
\label{xiid}
\end{eqnarray}
For giving the explicit form of $\Xi_{\io}^{(n-1)}$ and $\Xi_{\io}^{(n)}$,
we define the ``admissible patterns'' for the integer sequence
$\mu_1,\mu_2,\cd,\mu_k,\cd,\mu_n$ for $1 \leq k \leq n-1$ as follows:
\begin{eqnarray}
\begin{cases}
 1 \leq \mu_1 \leq n-1, \\
 0 \leq \mu_2 \leq \mu_1 -1, \\
\q\q\q  \cd, \\
 0 \leq \mu_k \leq \mu_{k-1}-1, \\
\q\q\q  \cd,
\end{cases}
\end{eqnarray}
where if $\mu_k$ does not exist, we define $\mu_k = 0$, and
\begin{eqnarray*}
S^{(\mu_1)}&:=& 
\begin{cases}
 id & (\mu_1 =1), \\
 S_{\mu_1-1;n-\mu_1} \cd S_{2;n-3}S_{1;n-2} & (\mu_1 \geq 2), \\
\end{cases} \\ 
S^{(\mu_k)}&:=&
\begin{cases}
  id & (\mu_k=0)\\
  S_{\mu_k+k-2;n-\mu_k} \cd S_{k+1;n-3}S_{k;n-2}S_{k-1;n} 
& ( \mu_k \geq 1\; \text{ and }\; k: \text{even }) \\
  S_{\mu_k+k-2;n-\mu_k} \cd S_{k+1;n-3}S_{k;n-2}S_{k-1;n-1} 
& ( \mu_k \geq 1\; \text{ and }\; k: \text{odd })
\end{cases}\; \text{ for }\; k \geq 2, \\
S^{(\mu_k')}&:=&
\begin{cases}
  id & (\mu_k=0)\\
  S_{\mu_k+k-2;n-\mu_k} \cd S_{k+1;n-3}S_{k;n-2}S_{k-1;n} 
& ( \mu_k \geq 1\; \text{ and }\; k: \text{odd }) \\
  S_{\mu_k+k-2;n-\mu_k} \cd S_{k+1;n-3}S_{k;n-2}S_{k-1;n-1} 
& ( \mu_k \geq 1\; \text{ and }\; k: \text{even })
\end{cases}\; \text{ for }\; k \geq 2, \\
\vp^{(\mu)}&:=& S^{(\mu_k)}\cd S^{(\mu_2)}S^{(\mu_1)}\;
\text{ for }\; \mu = (\mu_1,\mu_2,\cd,\mu_k,0,0,\cd), \\
\vp^{(\mu')}&:=& S^{(\mu_k')}\cd S^{(\mu_2')}S^{(\mu_1)}\;
\text{ for }\; \mu = (\mu_1,\mu_2',\cd,\mu_k',0,0,\cd).
\end{eqnarray*}
We prepare the symbol $X$ as follows:
\begin{eqnarray}
X &:=& x_{1;n-2} - x_{1;n-1}, \\
X' &:=& x_{1;n-2} - x_{1;n}.
\label{defXd}
\end{eqnarray}
For convenience, we define
\begin{eqnarray}
X_{j;i} &:=&
\begin{cases}
	x_{j;n} & \text{ if }\; j: \text{even and }\; i = n-1, \\
	x_{j;n-1}  & \text{ if }\; j: \text{even and }\; i = n, \\
        x_{j;i}  & \text{ other wise},
\end{cases} \\
X_{j;i}' &:=&
\begin{cases}
	x_{j;n} & \text{ if }\; j: \text{odd and }\; i = n-1, \\
	x_{j;n-1}  & \text{ if }\; j: \text{odd and }\; i = n, \\
        x_{j;i}  & \text{ other wise}.
\end{cases} \\
\end{eqnarray}
\begin{thm1} Let $\mu = (\mu_1,\mu_2,\cd,\mu_k,0,0,\cd)$ be the
admissible pattern.

\noindent
$(i)$ The forms $\vp^{(\mu)}X$ are given by
\begin{eqnarray}
\vp^{(\mu)}X = 
\begin{cases}
\sum_{k=1}^{l}
(X_{\mu_k+k-1;n-\mu_k-1} - X_{\mu_k+k-1;n-\mu_k})
& \text{ if }\; \mu_l = 1 , \\
\sum_{k=1}^{l}
(X_{\mu_k+k-1;n-\mu_k-1} - X_{\mu_k+k-1;n-\mu_k})
 + X_{l;n} & \text{ if }
\; \mu_{l} \geq 2.
\end{cases}
\label{sum_d}
\end{eqnarray}
where
$l := \text{max}\{ k : \mu_k \ne 0 \}$.
\vskip5pt

\noindent
$(ii)$ The forms $\vp^{(\mu')}X'$ are given by
\begin{eqnarray}
\vp^{(\mu')}X' = 
\begin{cases}
\sum_{k=1}^{l}
(X_{\mu_k+k-1;n-\mu_k-1}' - X_{\mu_k+k-1;n-\mu_k}')
& \text{ if }\; \mu_l' =1 , \\
\sum_{k=1}^{l}
(X_{\mu_k+k-1;n-\mu_k-1}' - X_{\mu_k+k-1;n-\mu_k}')
 + X_{l;n}' & \text{ if }\; \mu_l' \geq 2.
\end{cases}
\label{sum_dd}
\end{eqnarray}
where
$l := \text{max}\{ k : \mu_k' \ne 0 \}$.
\vskip5pt

\noindent
$($iii$)$
$\Xi_{\io}^{(n-1)}$ (resp. $\Xi_{\io}^{(n)}$) 
is the set of all linear forms which are consisted by
$\vp^{(\mu)}X$ (resp. $\vp^{(\mu')}X'$), 
where $\mu$ (resp. $\mu'$) are the admissible patterns.

\begin{proof}
We shall show (i), the case (ii) can be proved by the same argument.
First, we give a remark.
When $\mu_k - \mu_{k+1} =1$, the terms
$X_{\mu_k+k-1;n-\mu_k}$ and $X_{\mu_{k+1}+k;n-\mu_{k+1}-1}$
in the sum (\ref{sum_d}) are
canceled as $-X_{\mu_k+k-1;n-\mu_k} + X_{\mu_{k+1}+k;n-\mu_{k+1}-1}=0$.
We show the theorem by the induction on 
$|\mu| = \mu_1 + \mu_2 + \cd + \mu_i$
for $\mu = (\mu_1,\mu_2,\cd,\mu_i,0,0,\cd)$.

If $|\mu| = 1$, then $l = 1$ and
the sum of right hand side of (\ref{sum_d}) is
$X_{1;n-2}-X_{1;n-1} = x_{1;n-2}-x_{1;n-1}$ and equal to $\vp^{(\mu)}X$
by (\ref{defXd}).
We assume that $|\mu| = \mu_1 + \mu_2 + \cd + \mu_i = k-1$ for
$\mu = (\mu_1,\mu_2,\cd,\mu_i,0,0,\cd)$.
We consider the two cases: I) $\mu_i \rightarrow \mu_i +1$,
II) $``\mu_{i+1}=0'' \rightarrow ``\mu_{i+1}=1''$.
\vskip5pt

\noindent
I) $\mu_i \rightarrow \mu_i +1$.

\noindent
We shall prove the case $\mu_i=1$, the case $\mu_i \geq 2$ can be proved by
the same argument.
We have $l = i$ and
note that $1 \leq \mu_i \leq \mu_{i-1} -2$ by admissible
pattern of $\mu$ and the term $X_{\mu_i+i-1;n-\mu_i-1}$ in the sum
(\ref{sum_b}) is not canceled. We set 
$\mu' := (\mu_1,\mu_2,\cd,\mu_i+1,0,0,\cd)$.
We have
\begin{eqnarray*}
\vp^{(\mu')}X = S_{\mu_i+i-1;n-\mu_i-1}\vp^{(\mu)}X.
\end{eqnarray*}
The right hand side of (\ref{sum_d}) is
\begin{eqnarray*}
&&S_{\mu_i+i-1;n-\mu_i-1}(\sum_{k=1}^{l}
  (X_{\mu_k+k-1;n-\mu_k-1} - X_{\mu_k+k-1;n-\mu_k})) \\
=&& \sum_{k=1}^{l}
  (X_{\mu_k+k-1;n-\mu_k-1} - X_{\mu_k+k-1;n-\mu_k}) \\
 &&   - (X_{\mu_i+i-1;n-\mu_i-1} - X_{\mu_i+i-1;n-\mu_i}
     -X_{\mu_i+i;n-\mu_i-2} + X_{\mu_i+i;n-\mu_i-1}) \\
=&& \vp^{(\mu')}X.
\end{eqnarray*}

\noindent
II) $``\mu_{i+1}=0'' \rightarrow ``\mu_{i+1}=1''$.

\noindent
We set 
$\mu' := (\mu_1,\mu_2,\cd,\mu_i,1,0,\cd)$ $(i.e. \;\mu_{i+1}=1)$
and then $l = i+1$.
By the admissible pattern of $\mu$, we have $0 \leq \mu_{i+1}
\leq \mu_i -1$. This shows that $\mu_i \geq 2$ and
the term $X_{i;n}$ in the sum
(\ref{sum_d}) is not canceled.
If $i$ is odd, then 
we have
\begin{eqnarray*}
\vp^{(\mu')}X = S_{i;n}\vp^{(\mu)}X.
\end{eqnarray*}

The right hand side of (\ref{sum_d}) is
\begin{eqnarray*}
&&S_{i;n}(\sum_{k=1}^{l}
  (X_{\mu_k+k-1;n-\mu_k-1} - X_{\mu_k+k-1;n-\mu_k})+X_{i;n}) \\
=&& \sum_{k=1}^{l}
  (X_{\mu_k+k-1;n-\mu_k-1} - X_{\mu_k+k-1;n-\mu_k})+X_{i;n}
   - (X_{i;n} - X_{i+1;n-2} + X_{i+1;n}) \\
=&& \vp^{(\mu')}X.
\end{eqnarray*}
If $i$ is even, we can prove II) by the same argument.
\vskip5pt

\noindent
(iii) We sufficiently need to show that $\Xi_{\io}^{(n-1)}$ 
is closed under the
actions of all $S_{j;i}$, the other case can be proved by the
same argument.
We set for $1 \leq k \leq l+1$

$\mu =(\mu_1,\mu_2,\cd,\mu_k,\cd,\mu_l,0,0,\cd)$, 
$\mu^{+} =(\mu_1,\mu_2,\cd,\mu_k+1,\cd,\mu_l,0,0,\cd)$,

$\mu^{-} =(\mu_1,\mu_2,\cd,\mu_k-1,\cd,\mu_l,0,0,\cd)$.

\noindent
Then we have
\begin{eqnarray*}
S_{j;i}\vp^{(\mu)}X  =
\begin{cases}
\vp^{(\mu^{+})}X 
& \text{ if }\; (j;i)=(\mu_k+k-1;n-\mu_k) \text{ and }\;
\mu_{k-1}-\mu_k \ne 1, \\
\vp^{(\mu^{-})}X 
& \text{ if }\; (j;i)=(\mu_k+k-1;n-\mu_k+1) \text{ and }\;
\mu_{k}-\mu_{k+1} \ne 1, \\
\vp^{(\mu)}X & \text{ other wise. }
\end{cases}
\end{eqnarray*}
\end{proof}
\end{thm1}
Therefore, $\io$ satisfies the strict positivity assumption
by the explicit forms of $\Xi_{\io}^{(1,n-2)}$, $\Xi_{\io}^{(n-1)}$
and $\Xi_{\io}^{(n)}$. 
By the Lemma \ref{lmxi}, this shows that
\begin{eqnarray*}
\Xi_{\io}^{(1,n-2)}[\lm] &=&
 \{ \lm_i + x_{j;i-j} - x_{j;i-j+1} 
 : 1 \leq i \leq n-1, \; 1 \leq j \leq i \}, \\
\Xi_{\io}^{(k)}[\lm] &=&
 \{ \lm_k + \vp(\vec x) : \vp(\vec x) \in \Xi_{\io}^{(k)} :
 \mu \; \text{ are the admissible patterns }\} \; \text{ for } k=n-1,n
\end{eqnarray*}
by (\ref{xiid}), (\ref{sum_d}) and
$\Sigma_{\io}[\lm]$ is the polyhedral realization of $B(\lm)$.


\subsection{$F_4$ case}

We fix the $\io$ as follows:
$$\io := (\cd,\underbrace{4,3,2,1},\cd,
\underbrace{4,3,2,1},\underbrace{4,3,2,1}).$$
We define
\begin{eqnarray*}
\Xi_{\io} &:=& \{ S_{m_l}\cd S_{m_2}S_{m_1}(x_{j;1}) : 
l \geq 0, \; m_1,\cd,m_l \geq 1, \;1 \leq j \leq 6 \}, \\
\Sigma_{\io} &:=& \{\; \vec x \in \ZZ^{\ify}_{\io} \;|\;
\vp(\vec x) \geq 0 \text{ for any } \vp \in \Xi_{\io} \;\}.
\end{eqnarray*}

Using Theorem \ref{nz_modi}, we shall give the polyhedral realization
of $B(\ify)$. We give the explicit form of $\Xi_{\io}$ by direct
calculation and show that the polyhedral realization of $B(\ify)$
is equal to $\Sigma_{\io}$ along the following three steps:

\noindent
(i) We check the positivity assumption,

\noindent
(ii) For $\vec x = (\cd, x_{1;2},x_{1;1}) \in \Sigma_{\io}$,
we check that $x_{j;i} \geq 0$ for all $i$, $j$ $\geq 1$,

\noindent
((i), (ii) shows that $\Sigma_{\io}$ is
the polyhedral realization of $B(\ify)$),

\noindent
(iii) We determine when $x_{j;i} \equiv 0$ for convenience.

\begin{thm1}
We give the explicit form of the $\Xi_{\io}$ as follows:

\[
 \left\{
\begin{array}{ccc}
  \begin{array}{cccc}
 x_{j;1}  & x_{j;2} - x_{j+1;1} & 2x_{j;3} - x_{j+1;2}
& 2x_{j;4} - x_{j+3;1} \\
 x_{j+1;2} - 2x_{j+1;4} & x_{j+2;2} - x_{j+3;2}
& x_{j+2;1}- x_{j+4;1} & 2x_{j+2;4} - x_{j+4;2} \\
 x_{j+3;1} - x_{j+3;4} & x_{j+4;2} - 2x_{j+4;3}
& x_{j+5;1} - x_{j+5;2} & - x_{j+6;1}
  \end{array}
\\
  \begin{array}{ccc}
 2x_{j;4}+x_{j+1;2} - 2x_{j+1;3}
& 2x_{j;4} - x_{j+2;1} - x_{j+2;2}
& 2x_{j+1;3} - 2x_{j+1;4} - x_{j+3;1} \\
 x_{j+2;1} + x_{j+2;2} - 2x_{j+2;3}
& 2x_{j+2;2} - 2x_{j+2;3}-x_{j+3;1}
& x_{j+2;1}+x_{j+3;1} - x_{j+3;2} \\
 2x_{j+2;3} + x_{j+3;1} - 2x_{j+3;2}
& 2x_{j+2;4}+x_{j+3;1} - 2x_{j+3;3}
& 2x_{j+2;3}- x_{j+3;2} - x_{j+4;1} \\
x_{j+2;2} - x_{j+3;1}-x_{j+4;1}
& x_{j+3;2} - 2x_{j+3;4}-x_{j+4;1}
& 2x_{j+3;3} - 2x_{j+3;4}-x_{j+4;2}
  \end{array}
\\
  \begin{array}{cc}
 2x_{j+1;3}+x_{j+2;1} - 2x_{j+1;4}-x_{j+2;2}
& 2x_{j+2;4}+x_{j+3;2} - 2x_{j+3;3}-x_{j+4;1}
  \end{array}
\end{array}
 \right\}. \]
\end{thm1}

\noindent
(i) By the forms of $\Sigma_{\io}$, coefficients of $x_{1;1}$, $x_{1;2}$, 
$x_{1;3}$, $x_{1;4}$ are positive. This shows that ``positivity 
assumption'' is satisfied.
\vskip5pt

\noindent
(ii)
We assume $x_{j;1} \geq 0$ for any $j \geq 1$.
We show that $x_{j;2}$, $x_{j;3}$, $x_{j;4} \geq 0$ for any $j \geq 1$.
By the results of $(i)$, we have
$$  x_{j;2} \geq x_{j+1;1}, \;\; 2x_{j;3} \geq x_{j+1;2},
\;\; 2x_{j;4} \geq x_{j+3;1}.$$
Then we have $x_{j;2} \geq 0$
for $j \geq 1$ since $x_{j;2} \geq x_{j+1;1}$. 
Similarly, we obtain $x_{j;3} \geq 0$,
$x_{j;4} \geq 0$ for $j \geq 1$ since
$2x_{j;3} \geq x_{j+1;2}$,
$2x_{j;4} \geq x_{j+3;1}$ respectively.
%
%
%
%
%
These $(i)$, $(ii)$, $(iii)$ shows that $\Sigma_{\io}$ is the
polyhedral realization of $B(\ify)$.
\vskip5pt

\noindent
(iii) We determine when $x_{j;i} \equiv 0$.
By the forms of $\Sigma_{\io}$ and (ii), we have
$$0 \geq x_{j+6;1} \geq 0.$$
This shows $x_{m;1} \equiv 0$ for $m \geq 7$.
Similarly, we obtain
$$x_{m;2} \equiv 0, \;\;x_{m;3} \equiv 0,\;\,
x_{m;4} \equiv 0 \;\text{ for } m \geq 7$$
since $x_{j+5;1} \geq x_{j+5;2}$, $x_{j+4;1} \geq x_{j+4;3}$,
$x_{j+3;1} \geq x_{j+3;4}$ respectively.
In particular, the parameter $j$ in $\Xi_{\io}$ is $1 \leq j \leq 6$.



Next, we give the polyhedral realization of $B(\lm)$ for
$\lm := \lm_1\Lambda_1 + \cd + \lm_4\Lambda_4$.
Here we set for $1 \leq i \leq 4$
\begin{eqnarray*}
\Xi_{\io}^{(i)} &:=& 
\{ S_{j_k}\cd S_{j_2}  S_{j_1}\xi^{(i)}
(\vec x) : k \geq 0, \;j_1,\cd,j_k \geq 1 \},\\
\Xi_{\io}^{(i)}[\lm] &:=& 
\{ \hat S_{j_k}\cd \hat S_{j_2} \hat S_{j_1}\lm^{(i)}
(\vec x) : k \geq 0, \; j_1,\cd,j_k \geq 1 \},\\
\Xi_{\io}[\lm] &:=& \Xi_{\io} \cup \Xi_{\io}^{(1)}[\lm]\cup \cd
\cup\Xi_{\io}^{(4)}[\lm] \\
\Sigma_{\io}[\lm] &:=& \{ \vec x \in \ZZ_{\io}[\lm] : \vp(\vec x)
\geq 0\; \text{ for any }\; \vp \in \Xi_{\io}[\lm] \}.
\end{eqnarray*}
In order to show that $\Sigma_{\io}[\lm]$ is the polyhedral
realization of $B(\lm)$, we give the
explicit forms of $\Xi_{\io}^{(i)}$ for $1 \leq i \leq 4$
by direct calculation as follws:
\begin{eqnarray*}
\Xi_{\io}^{(1)}&=& \{ -x_{1;1} \}, \\
\Xi_{\io}^{(2)}&=& \{ x_{1;1}-x_{1;2}, \; -x_{2;1} \}, \\
\Xi_{\io}^{(3)}&=& \\
&& \left\{
\begin{array}{ccc}
  \begin{array}{cccc}
 x_{1;2} - x_{1;3} & x_{1;3} - x_{3;1}
& x_{2;1} - x_{2;4} & x_{2;3} - x_{3;3} \\
 x_{1;4}- x_{3;4} & x_{2;4} - x_{5;1}
& x_{4;1} - x_{4;3} & x_{4;3} - x_{5;2} \\
 x_{4;4} - x_{5;3} & - x_{5;4}
&&
  \end{array}
\\
  \begin{array}{ccc}
 x_{1;3}+x_{2;1} - x_{2;2}
& x_{1;4} + x_{2;1} - x_{2;3}
& x_{1;4} + x_{2;3} - x_{3;2} \\
 x_{1;4} + x_{2;4} - x_{3;3}
& x_{2;4} + x_{3;2}-2x_{3;3}
& x_{2;4}+x_{4;1} - x_{4;2}
  \end{array}
\\
  \begin{array}{ccc}
 x_{2;2}-x_{2;4} - x_{3;1}
& 2x_{2;3} - x_{2;4} - x_{3;2}
& x_{3;2} - x_{3;3} - x_{3;4} \\
 x_{2;3} - x_{2;4} - x_{3;4}
& x_{3;3} - x_{3;4}-x_{5;1}
& x_{4;2}-x_{4;3} - x_{5;1}
  \end{array}
\\
  \begin{array}{cc}
 x_{1;4}+x_{2;2} - x_{2;3}-x_{3;1}
& x_{3;3}+x_{4;1} - x_{3;4}-x_{4;2}
  \end{array}
\end{array}
 \right\}, \\
\Xi_{\io}^{(4)}&=& \\
&& \left\{
\begin{array}{ccc}
  \begin{array}{cccc}
 x_{1;3} - x_{1;4} & x_{2;2} - x_{2;3}
& x_{2;3} - x_{4;1} & x_{3;1} - x_{3;4} \\
 x_{3;3}- x_{4;3} & x_{2;4} - x_{4;4}
& x_{3;4} - x_{6;1} & x_{5;1} - x_{5;3} \\
 x_{5;3} - x_{6;2} & x_{5;4} - x_{6;3}
& - x_{6;4} &
  \end{array}
\\
  \begin{array}{ccc}
 x_{2;3}+x_{3;1} - x_{3;2}
& x_{2;4} + x_{3;1} - x_{3;3}
& x_{2;4} + x_{3;3} - x_{4;2} \\
 x_{2;4} + x_{3;4} - x_{4;3}
& x_{3;4} + x_{4;2}-2x_{4;3}
& x_{3;4}+x_{5;1} - x_{5;2}
  \end{array}
\\
  \begin{array}{ccc}
 x_{3;2}-x_{3;4} - x_{4;1}
& 2x_{3;3} - x_{3;4} - x_{4;2}
& x_{4;2} - x_{4;3} - x_{4;4} \\
 x_{3;3} - x_{3;4} - x_{4;4}
& x_{4;3} - x_{4;4}-x_{6;1}
& x_{5;2}-x_{5;3} - x_{6;1}
  \end{array}
\\
  \begin{array}{cc}
 x_{2;4}+x_{3;2} - x_{3;3}-x_{4;1}
& x_{4;3}+x_{5;1} - x_{4;4}-x_{5;2}
  \end{array}
\end{array}
 \right\}. 
\end{eqnarray*}
Therefore, $\io$ satisfies the strict positivity assumption
by the explicit forms of $\Xi_{\io}^{(i)}$
for $1 \leq i \leq 4$. 
This shows that
\begin{eqnarray*}
\Xi_{\io}^{(1)}[\lm] &=&
 \{ \lm_1 - x_{1;1} \}, \\
\Xi_{\io}^{(2)}[\lm] &=&
 \{ \lm_2 + x_{1;1} - x_{1;2}, \; \lm_2 -x_{2;1} \}, \\
\Xi_{\io}^{(k)}[\lm] &=&
 \{ \lm_k + \vp_k(\vec x) :
 \vp_k(\vec x) \in \Xi_{\io}^{(k)}\}\; \text{ for } k=3,4.
\end{eqnarray*}
by the Lemma \ref{lmxi} and
$\Sigma_{\io}[\lm]$ is the polyhedral realization of $B(\lm)$.


\subsection{$E_6$ case}

We fix the $\io$ as follows:
$$\io := (\cd,\underbrace{6,5,4,3,2,1},\cd,
\underbrace{6,5,4,3,2,1},\underbrace{6,5,4,3,2,1}).$$
We define
\begin{eqnarray*}
\Xi_{\io} &:=& \{ S_{m_l}\cd S_{m_2}S_{m_1}(x_{j;1}) : 
l \geq 0,\; m_1,\cd,m_l \geq 1,\;1 \leq j \leq 8 \}, \\
\Sigma_{\io} &:=& \{\; \vec x \in \ZZ^{\ify}_{\io} :
\vp(\vec x) \geq 0 \text{ for any } \vp \in \Xi_{\io} \;\}.
\end{eqnarray*}

We give the explicit form of $\Xi_{\io}$ by direct calculation
and show that the polyhedral realization of $B(\ify)$ is equal to
$\Sigma_{\io}$
along the following steps
similar to the case of $F_4$:

\noindent
(i) We check the positivity assumption

\noindent
(ii) For $\vec x = (\cd, x_{1;2},x_{1;1}) \in \Sigma_{\io}$,
we check that $x_{j;i} \geq 0$ for all $i$, $j$ $\geq 1$

\noindent
(iii) We determine when $x_{j;i} \equiv 0$ for convenience.

\begin{thm1}
We give the explicit form of the
$\Xi_{\io}$ as follows:

\[
 \left\{
\begin{array}{ccc}
  \begin{array}{cccc}
 x_{j;1}  & x_{j;2} - x_{j+1;1} & x_{j;3} - x_{j+1;2}
& x_{j;4} - x_{j+1;6} \\
 x_{j;6} - x_{j+1;5} & x_{j+2;2} - x_{j+2;4}
& x_{j;5} - x_{j+4;1} & x_{j+3;1} - x_{j+3;6} \\
 x_{j+2;6} -x_{j+2;2} & x_{j+3;4}- x_{j+4;3}
& x_{j+3;5} - x_{j+4;4} &- x_{j+4;5}
  \end{array}
\\
  \begin{array}{ccc}
 x_{j;4}+x_{j;6} - x_{j+1;3}
& x_{j;5} + x_{j;6}- x_{j+1;4}
& x_{j;5} + x_{j+2;2} - x_{j+2;3} \\
 x_{j;5} + x_{j+3;1} - x_{j+3;2}
& x_{j+2;6} +x_{j+3;1} - x_{j+3;3}
& x_{j+1;3} - x_{j+1;5} -x_{j+1;6}\\
 x_{j+1;4}  - x_{j+1;5} - x_{j+4;1}
& x_{j+2;3} - x_{j+2;4} -x_{j+4;1}
& x_{j+3;2} - x_{j+3;6} - x_{j+4;1}
  \end{array}
\\
  \begin{array}{cc}
 x_{j+3;3} - x_{j+3;6}-x_{j+4;2}
&  x_{j;5}+x_{j+1;3}- x_{j+1;4}-x_{j+1;6}
  \end{array}
\\
  \begin{array}{cc}
 x_{j+1;4}+x_{j+2;2} - x_{j+1;5}-x_{j+2;3}
& x_{j+1;4}+x_{j+3;1} - x_{j+1;5}-x_{j+3;2} \\
 x_{j+2;3}+x_{j+3;1} - x_{j+2;4}-x_{j+3;2}
& x_{j+2;6}+x_{j+3;2} - x_{j+3;3}-x_{j+4;1}
  \end{array}
\end{array}
 \right\}. \]
\end{thm1}

\noindent
(i) By the form of $\Sigma_{\io}$, coefficients of $x_{1;1}$, $x_{1;2}$, 
$x_{1;3}$, $x_{1;4}$, $x_{1;5}$, $x_{1;6}$
are positive. This shows that ``positivity 
assumption'' is satisfied.
\vskip5pt

\noindent
(ii)
We assume $x_{j;1} \geq 0$ for any $j \geq 1$.
We show that $x_{j;2}$, $x_{j;3}$, $x_{j;4}$,
$x_{j;5}$, $x_{j;6} \geq 0$ for any $j \geq 1$.
By the form of $\Sigma_{\io}$, we have
$$x_{j;2} \geq x_{j+1;1}, \;\;x_{j;3} \geq x_{j+1;2},\;\;
x_{j;5} \geq x_{j+4;1},\;\;x_{j;6} \geq x_{j+1;5},\;\;
x_{j;4} \geq x_{j+1;6}.$$
This shows $x_{j;2} \geq 0$ for $ j\geq 1$ since $x_{j+1;1} \geq 0$ and
similarly, we have $x_{j;3} \geq 0$, $x_{j;5} \geq 0$,
$x_{j;6} \geq 0$, $x_{j;4} \geq 0$ for $j \geq 1$
since $x_{j;3} \geq x_{j+1;2}$,
$x_{j;5} \geq x_{j+4;1}$, $x_{j;6} \geq x_{j+1;5}$,
$x_{j;4} \geq x_{j+1;6}$ respectively.




%





%



\vskip5pt

\noindent
(iii)
We determine when $x_{j;i} \equiv 0$.
We have
$$0 \geq x_{j+4;5} \geq 0.$$
This shows $x_{m;5} \equiv 0$ for $m \geq 5$.
Similarly, we have
$x_{m;4} \equiv 0$ for $m \geq 6$,
$x_{m;3} \equiv 0$ for $m \geq 7$,
$x_{m;6} \equiv 0$ for $m \geq 7$,
$x_{m;2} \equiv 0$ for $m \geq 8$,
$x_{m;1} \equiv 0$ for $m \geq 9$
since $x_{j+3;4} \geq x_{j+4;3}$,
$x_{j+3;5} \geq x_{j+4;4}$,
$x_{j;4} \geq x_{j+1;6}$,
$x_{j;3} \geq x_{j+1;2}$,
$x_{j;2} \geq x_{j+1;1}$
respectively.
In particular, the parameter $j$ in $\Xi_{\io}$ runs
$1 \leq j \leq 8$.



Next, we give the polyhedral realization of $B(\lm)$ for
$\lm := \lm_1\Lambda_1 + \cd + \lm_6\Lambda_6$.
Here we set for $1 \leq i \leq 6$
\begin{eqnarray*}
\Xi_{\io}^{(i)} &:=& 
\{ S_{j_k}\cd S_{j_2}  S_{j_1}\xi^{(i)}
(\vec x) : k \geq 0, \;j_1,\cd,j_k \geq 1 \},\\
\Xi_{\io}^{(i)}[\lm] &:=& 
\{ \hat S_{j_k}\cd \hat S_{j_2} \hat S_{j_1}\lm^{(i)}
(\vec x) : k \geq 0, \; j_1,\cd,j_k \geq 1 \},\\
\Xi_{\io}[\lm] &:=& \Xi_{\io} \cup \Xi_{\io}^{(1)}[\lm]\cup \cd
\cup\Xi_{\io}^{(6)}[\lm] \\
\Sigma_{\io}[\lm] &:=& \{ \vec x \in \ZZ_{\io}[\lm] : \vp(\vec x)
\geq 0\; \text{ for any }\; \vp \in \Xi_{\io}[\lm] \}.
\end{eqnarray*}
In order to show that $\Sigma_{\io}[\lm]$ is the polyhedral
realization of $B(\lm)$, we give the
explicit forms of $\Xi_{\io}^{(i)}$ for $1 \leq i \leq 6$
by direct calculation as follws:
\begin{eqnarray*}
\Xi_{\io}^{(1)}&&= \{ -x_{1;1} \}, \\
\Xi_{\io}^{(2)}&&= \{ x_{1;1}-x_{1;2}, \; -x_{2;1} \}, \\
\Xi_{\io}^{(3)}&&= \{ x_{1;2}-x_{1;3}, \; x_{2;1} - x_{2;2},\;
-x_{3;1} \}, \\
\Xi_{\io}^{(4)}&&= \\
&& \left\{
\begin{array}{ccc}
  \begin{array}{ccccccc}
 x_{1;3} - x_{1;4} & x_{2;2} - x_{2;6}
& x_{1;6} - x_{4;1} & x_{3;1} - x_{3;5}
& x_{2;4}- x_{4;2} & x_{2;5} - x_{4;6}
& x_{3;6} - x_{4;4} \\
 x_{5;2} - x_{5;3} & x_{6;1} - x_{6;2} 
& -x_{7;1} &&&&
  \end{array}
\\
  \begin{array}{cccc}
 x_{1;6}+x_{2;2} - x_{2;3}
& x_{1;6} + x_{3;1} - x_{3;2}
& x_{2;4} + x_{3;1} - x_{3;3}
& x_{2;5} + x_{3;1} - x_{3;4} \\
 x_{2;5} + x_{3;6}-x_{4;3} &&&
  \end{array}
\\
  \begin{array}{cccc}
 x_{2;3}-x_{2;6} - x_{4;1}
& x_{3;2} - x_{3;5} - x_{4;1}
& x_{3;3} - x_{3;5} - x_{4;2}
& x_{3;4} - x_{3;5} - x_{4;6} \\
 x_{4;3} - x_{4;4}-x_{4;6}
&&&
  \end{array}
\\
  \begin{array}{ccc}
 x_{2;3}+x_{3;1} - x_{2;6}-x_{3;2}
& x_{2;4}+x_{3;2} - x_{3;3}-x_{4;1}
& x_{2;5}+x_{3;2} - x_{3;4}-x_{4;1} \\
 x_{2;5}+x_{3;3} - x_{3;4}-x_{4;2}
& x_{3;4}+x_{3;6} - x_{3;5}-x_{4;3}
&
  \end{array}
\end{array}
 \right\}, \\
\Xi_{\io}^{(5)}&&= \\
&& \left\{
\begin{array}{ccc}
  \begin{array}{ccccccc}
 x_{1;4} - x_{1;5} & x_{2;3} - x_{2;4} 
& x_{3;2} - x_{3;6} & x_{2;6} - x_{5;1} 
& x_{4;1} - x_{4;5} & x_{3;4}- x_{5;2} 
& x_{3;5} - x_{5;6} \\
 x_{4;6} - x_{5;4} & x_{6;2} - x_{6;3} 
& x_{7;1} - x_{7;2} & -x_{8;1} &&&
  \end{array}
\\
  \begin{array}{cccc}
 x_{2;6}+x_{3;2} - x_{3;3}
& x_{2;6} + x_{4;1} - x_{4;2}
& x_{3;4} + x_{4;1} - x_{4;3}
& x_{3;5} + x_{4;1} - x_{4;4} \\
 x_{3;5} + x_{4;6}-x_{5;3}
&&&
  \end{array}
\\
  \begin{array}{cccc}
 x_{3;3}-x_{3;6} - x_{5;1}
& x_{4;2} - x_{4;5} - x_{5;1}
& x_{4;3} - x_{4;5} - x_{5;2}
& x_{4;4} - x_{4;5} - x_{5;6} \\
 x_{5;3} - x_{5;4}-x_{5;6} 
&&&
  \end{array}
\\
  \begin{array}{ccc}
 x_{3;3}+x_{4;1} - x_{3;6}-x_{4;2}
& x_{3;4}+x_{4;2} - x_{4;3}-x_{5;1}
& x_{3;5}+x_{4;2} - x_{4;4}-x_{5;1} \\
 x_{3;5}+x_{4;3} - x_{4;4}-x_{5;2}
& x_{4;4}+x_{4;6} - x_{4;5}-x_{5;3}
&
  \end{array}
\end{array}
 \right\}, \\
\Xi_{\io}^{(6)}&&= \\
&& \left\{
\begin{array}{ccc}
  \begin{array}{ccccccc}
 x_{1;3} - x_{1;6} & x_{2;2} - x_{2;5} 
& x_{1;4} - x_{4;1} 

& x_{3;3} - x_{4;3} 
& x_{3;1} - x_{6;1} & x_{2;6}- x_{4;6} 
& x_{1;5} - x_{4;5} \\

 x_{5;1} - x_{5;4} & x_{3;5} - x_{6;2} 
& x_{5;6} - x_{6;3} & -x_{6;1} &&&
  \end{array}
\\
  \begin{array}{cccc}
 x_{1;4}+x_{2;2} - x_{2;3}
& x_{1;5} + x_{2;2} - x_{2;4}
& x_{1;4} + x_{3;1} - x_{3;2}
& x_{1;5} + x_{3;1} - x_{3;6} \\
 x_{2;6}+x_{3;1} - x_{2;4}
& x_{1;5} + x_{2;6} - x_{4;2}

& x_{2;1} + x_{4;2} - x_{4;3}
& x_{1;5} + x_{3;4} - x_{4;3} \\
 x_{3;1}+x_{5;1} - x_{5;2}
& x_{2;6} + x_{3;6} - x_{4;3}
& x_{1;5} + x_{3;5} - x_{4;4}

& x_{3;5} + x_{5;1} - x_{5;3}
  \end{array}
\\
  \begin{array}{cccc}
 x_{2;3}-x_{2;5} - x_{4;1}

& x_{3;2} - x_{4;1} - x_{6;1}
& x_{3;3} - x_{3;6} - x_{4;6}
& x_{2;4} - x_{2;5} - x_{4;5} \\
 x_{3;3}-x_{4;2} - x_{6;1}
& x_{3;3} - x_{3;4} - x_{4;5}
& x_{4;2} - x_{4;5} - x_{4;6}
& x_{3;4} - x_{4;6} - x_{6;1} \\
 x_{3;6}-x_{4;5} - x_{6;1}

& x_{5;2} - x_{5;4} - x_{6;1}
& x_{4;4} - x_{4;5} - x_{6;2}
& x_{5;3} - x_{5;4} - x_{6;2} \\
  \end{array}
\\
  \begin{array}{ccc}
 x_{2;3}+x_{3;1} - x_{2;5}-x_{3;2}
& x_{1;5}+x_{2;3} - x_{2;4}-x_{4;1}
& x_{2;4}+x_{3;1} - x_{2;5}-x_{3;6} \\
 x_{1;5}+x_{3;2} - x_{3;6}-x_{4;1}
& x_{3;1}+x_{3;3} - x_{3;4}-x_{3;6}
& x_{2;6}+x_{3;2} - x_{3;4}-x_{4;1} \\
 x_{2;4}+x_{2;5} - x_{2;6}-x_{4;2}

& x_{1;5}+x_{3;3} - x_{3;6}-x_{4;2}
& x_{2;6}+x_{3;3} - x_{3;4}-x_{4;2} \\
 x_{3;2}+x_{4;2} - x_{4;1}-x_{4;3}
& x_{2;4}+x_{3;4} - x_{2;5}-x_{4;3}
& x_{3;2}+x_{5;1} - x_{4;1}-x_{5;2} \\
 x_{2;4}+x_{3;5} - x_{2;5}-x_{4;4}
& x_{3;4}+x_{4;2} - x_{4;3}-x_{4;6}
& x_{3;3}+x_{4;2} - x_{5;1}-x_{5;2} \\

 x_{3;3}+x_{3;5} - x_{3;4}-x_{4;4}
& x_{3;6}+x_{4;2} - x_{4;3}-x_{4;5}
& x_{3;5}+x_{4;2} - x_{4;4}-x_{4;6} \\
 x_{3;4}+x_{5;1} - x_{4;6}-x_{5;2}
& x_{3;4}+x_{3;6} - x_{4;3}-x_{6;1}
& x_{3;6}+x_{5;1} - x_{4;5}-x_{5;2} \\
 x_{3;5}+x_{3;6} - x_{4;4}-x_{6;1}
& x_{4;4}+x_{5;1} - x_{4;5}-x_{5;3}

& x_{3;5}+x_{5;2} - x_{5;3}-x_{6;1}
  \end{array}
\\
  \begin{array}{ccc}
 x_{1;5}+x_{2;6} + x_{3;1}-x_{3;3}
& x_{3;4}+x_{3;6} + x_{4;2}-2x_{4;3}

& 2x_{3;3}-x_{3;4} - x_{3;6}-x_{4;2} \\
 x_{4;3}-x_{4;5} - x_{4;6}-x_{6;1}
  \end{array}
\\
  \begin{array}{cc}
 x_{2;4}+x_{2;6} + x_{3;2}-x_{2;5}
 -x_{3;3}-x_{4;1} 
& x_{3;5}+x_{4;3}+ x_{5;1}-x_{4;4}
 -x_{4;6}-x_{5;2}
  \end{array}
\\
  \begin{array}{cc}
  x_{1;5} +x_{2;3}+x_{3;1} -x_{2;4} - x_{3;2}
& x_{2;4} +x_{2;6}+x_{3;1} -x_{2;5}- x_{3;3} \\
  x_{1;5} +x_{2;6}+x_{3;2} -x_{3;3} - x_{4;1}

&  x_{3;5} +x_{3;6}+x_{4;2} -x_{4;3} - x_{4;4} \\
  x_{3;4} +x_{3;6}+x_{5;1} -x_{4;3}- x_{5;2}

&  x_{3;5} +x_{3;6}+x_{5;1} -x_{4;4} - x_{5;2} \\
  x_{2;4} +x_{3;2}-x_{2;5} -x_{3;6} - x_{4;1}

& x_{3;2} +x_{3;3}-x_{3;4} -x_{3;6}- x_{4;1} \\
  x_{2;4} +x_{3;3}-x_{2;5} -x_{3;6} - x_{4;2}

&  x_{4;3} +x_{5;1}-x_{4;5} -x_{4;6} - x_{5;2} \\
 x_{3;5} +x_{4;3}-x_{4;4} -x_{4;6}- x_{6;1}

& x_{4;4} +x_{5;2}-x_{4;5} -x_{5;3} - x_{6;1}
  \end{array}
\end{array}
 \right\}.
\end{eqnarray*}
Therefore, $\io$ satisfies the strict positivity assumption
by the explicit forms of $\Xi_{\io}^{(i)}$
for $1 \leq i \leq 6$. 
This shows that
\begin{eqnarray*}
\Xi_{\io}^{(1)}[\lm] &=&
 \{ \lm_1 - x_{1;1} \}, \\
\Xi_{\io}^{(2)}[\lm] &=&
 \{ \lm_2 + x_{1;1} - x_{1;2}, \; \lm_2 -x_{2;1} \}, \\
\Xi_{\io}^{(3)}[\lm] &=&
 \{ \lm_3 + x_{1;2} - x_{1;3}, \; \lm_2 +x_{2;1}-x_{2;2},
\; \lm_3 - x_{3;1} \}, \\
\Xi_{\io}^{(k)}[\lm] &=&
 \{ \lm_k + \vp_k(\vec x) :
 \vp_k(\vec x) \in \Xi_{\io}^{(k)}\}\; \text{ for } k=4,5,6.
\end{eqnarray*}
by the Lemma \ref{lmxi} and
$\Sigma_{\io}[\lm]$ is the polyhedral realization of $B(\lm)$.


\subsection{$E_7$ case}

We fix the $\io$ as follows:
$$\io := (\cd,\underbrace{7,6,5,4,3,2,1},\cd,
\underbrace{7,6,5,4,3,2,1},\underbrace{7,6,5,4,3,2,1}).$$
We define
\begin{eqnarray*}
\Xi_{\io} &:=& \{ S_{m_l}\cd S_{m_2}S_{m_1}(x_{j;1}) : 
l \geq 0,\; m_1,\cd,m_l \geq 1,\;1 \leq j \leq 9 \}, \\
\Sigma_{\io} &:=& \{\; \vec x \in \ZZ^{\ify}_{\io} :
\vp(\vec x) \geq 0 \text{ for any } \vp \in \Xi_{\io} \;\}.
\end{eqnarray*}

\begin{thm1}
We give the explicit form of the
$\Xi_{\io}$ as follows:

\[
 \left\{
\begin{array}{ccc}
  \begin{array}{cccc}
 x_{j;1}  & x_{j;2} - x_{j+1;1} & x_{j;3} - x_{j+1;2}
& x_{j;4} - x_{j+1;3} \\
 x_{j;5} - x_{j+1;7} & x_{j;7} - x_{j+1;6}
& x_{j+2;3} - x_{j+2;5} & x_{j;6} - x_{j+5;1} \\
 x_{j+3;2} - x_{j+3;7} & x_{j+2;7} - x_{j+5;2}
& x_{j+4;1} - x_{j+4;6} &  x_{j+3;5} - x_{j+5;3} \\
 x_{j+3;6} - x_{j+5;7} & x_{j+4;7} - x_{j+5;5}
& x_{j+6;3} - x_{j+6;4} & x_{j+7;2} - x_{j+7;3} \\
 x_{j+8;1} - x_{j+8;2} &  - x_{j+9;1}
& &
  \end{array}
\\
  \begin{array}{ccc}
 x_{j;5}+x_{j;7} - x_{j+1;4}
& x_{j;6} + x_{j;7}-  x_{j+1;5}
& x_{j;6} + x_{j+2;3} - x_{j+2;4} \\
 x_{j;6} + x_{j+3;2} - x_{j+3;3}
& x_{j;6} +x_{j+4;1} - x_{j+4;2}
& x_{j+2;7} +x_{j+3;2} - x_{j+3;4} \\
 x_{j+2;7} + x_{j+4;1} - x_{j+4;3}
& x_{j+3;5} + x_{j+4;1} - x_{j+4;4}
& x_{j+3;6} +x_{j+4;1} - x_{j+4;5} \\
 x_{j+3;6} +x_{j+4;7} - x_{j+5;4}
& x_{j+1;4}  - x_{j+1;6} -x_{j+1;7}
& x_{j+1;5}  - x_{j+1;6} -x_{j+5;1}\\
 x_{j+2;4}  - x_{j+2;5} - x_{j+5;1}
& x_{j+3;3} - x_{j+3;7} -x_{j+5;1}
& x_{j+3;4} - x_{j+3;7} - x_{j+5;2}\\
 x_{j+4;2}  - x_{j+4;6} - x_{j+5;1}
& x_{j+4;3} - x_{j+4;6} -x_{j+5;2}
& x_{j+4;4} - x_{j+4;6}  -x_{j+5;3}\\
 x_{j+4;5}  - x_{j+4;6} - x_{j+5;7}
& x_{j+5;4} - x_{j+5;5} -x_{j+5;7}
&
  \end{array}
\\
  \begin{array}{cc}
 x_{j;6}+x_{j+1;4} - x_{j+1;5}-x_{j+1;7}
& x_{j+1;5}+x_{j+2;3} - x_{j+1;6}-x_{j+2;4} \\
 x_{j+1;5}+x_{j+3;2} - x_{j+1;6}-x_{j+3;3}
& x_{j+2;4}+x_{j+3;2} -x_{j+2;5}-x_{j+3;3}\\
 x_{j+1;5}+x_{j+4;1} - x_{j+1;6}-x_{j+4;2}
& x_{j+2;4}+x_{j+4;1} - x_{j+2;5}-x_{j+4;2} \\
 x_{j+3;3}+x_{j+4;1} - x_{j+3;7}-x_{j+4;2}
& x_{j+2;7}+x_{j+3;3} -x_{j+3;4}-x_{j+5;1}\\
 x_{j+3;4}+x_{j+4;1} - x_{j+3;7}-x_{j+4;3}
& x_{j+2;7}+x_{j+4;2} - x_{j+4;3}-x_{j+5;1} \\
 x_{j+3;5}+x_{j+4;2} - x_{j+4;4}-x_{j+5;1}
& x_{j+3;6}+x_{j+4;2} - x_{j+4;5}-x_{j+5;1}\\
 x_{j+3;5}+x_{j+4;3} - x_{j+4;4}-x_{j+5;2}
& x_{j+3;6}+x_{j+4;3} - x_{j+4;5}-x_{j+5;2} \\
 x_{j+3;6}+x_{j+4;4} -x_{j+4;5}-x_{j+5;3}
& x_{j+4;5}+x_{j+4;7} - x_{j+4;6}-x_{j+5;4}\\
 x_{j+2;7}+x_{j+3;3}+x_{j+4;1} - x_{j+3;4}-x_{j+4;2}
& x_{j+3;4}+x_{j+4;2} - x_{j+3;7}-x_{j+4;3}-x_{j+5;1}
  \end{array}
\end{array}
 \right\}. \]
\end{thm1}

\noindent
(i) By the form of $\Sigma_{\io}$, coefficients of $x_{1;1}$, $x_{1;2}$, 
$x_{1;3}$, $x_{1;4}$, $x_{1;5}$, $x_{1;6}$, $x_{1;7}$
are positive. This shows that ``positivity 
assumption'' is satisfied.
\vskip5pt

\noindent
$(ii)$
We assume $x_{j;1} \geq 0$ for any $j \geq 1$.
We show that $x_{j;2}$, $x_{j;3}$, $x_{j;4}$,
$x_{j;5}$, $x_{j;6}$, $x_{1;7} \geq 0$ for any $j \geq 1$.
By the form of $\Sigma_{\io}$, we have
$$x_{j;2} \geq x_{j+1;1}, \;\;x_{j;3} \geq x_{j+1;2},\;\;
x_{j;4} \geq x_{j+1;3},\;\;
x_{j;6} \geq x_{j+5;1},\;\;
x_{j;7} \geq x_{j+1;6},\;\;
x_{j;5} \geq x_{j+1;7}.
$$
This shows $x_{j;2} \geq 0$ for $ j\geq 1$ since $x_{j+1;1} \geq 0$ and
similarly, we have $x_{j;3} \geq 0$, $x_{j;4} \geq 0$,
$x_{j;6} \geq 0$, $x_{j;7} \geq 0$, $x_{j;5} \geq 0$
for $j \geq 1$
since $x_{j;3} \geq x_{j+1;2}$,
$x_{j;4} \geq x_{j+1;3}$, $x_{j;6} \geq x_{j+5;1}$,
$x_{j;7} \geq x_{j+1;6}$, $x_{j;5} \geq x_{j+1;7}$ respectively.
\vskip5pt













\noindent
$(iii)$
We determine when $x_{j;i} \equiv 0$.
We have
$$0 \geq x_{j+9;1} \geq 0.$$
This shows $x_{m;1} \equiv 0$ for $m \geq 10$.
Similarly, we have
$x_{m;2} \equiv 0$, 
$x_{m;3} \equiv 0$,
$x_{m;4} \equiv 0$,
$x_{m;5} \equiv 0$,
$x_{m;6} \equiv 0$,
$x_{m;7} \equiv 0$ for $m \geq 10$
since $x_{j+8;1} \geq x_{j+8;2}$,
$x_{j+7;2} \geq x_{j+7;3}$,
$x_{j+6;3} \geq x_{j+6;4}$,
$x_{j+2;3} \geq x_{j+2;5}$,
$x_{j+4;1} \geq x_{j+4;6}$,
$x_{j+3;2} \geq x_{j+3;7}$
respectively. In particular, the parameter $j$ in $\Xi_{\io}$ runs
$1 \leq j \leq 9$.


\subsection{$E_8$ case}

We fix the $\io$ as follows:
$$\io := (\cd,\underbrace{8,7,6,5,4,3,2,1},\cd,
\underbrace{8,7,6,5,4,3,2,1},\underbrace{8,7,6,5,4,3,2,1}).$$
We define
\begin{eqnarray*}
\Xi_{\io} &:=& \{ S_{m_l}\cd S_{m_2}S_{m_1}(x_{j;1}) : 
l \geq 0,\; m_1,\cd,m_l \geq 1,\;1 \leq j \leq 15 \}, \\
\Sigma_{\io} &:=& \{\; \vec x \in \ZZ^{\ify}_{\io} :
\vp(\vec x) \geq 0 \text{ for any } \vp \in \Xi_{\io} \;\}.
\end{eqnarray*}

\begin{thm1}
We give the explicit form of the
$\Xi_{\io}$ as follows:

\[
\begin{array}{cccc}
  \begin{array}{ccccc}
 x_{j;1} & x_{j;2} - x_{j+1;1} & x_{j;3} - x_{j+1;2}
& x_{j;4} - x_{j+1;3} & x_{j;5} - x_{j+1;4} \\
 x_{j;6} - x_{j+1;8} & x_{j;7} - x_{j+6;1} & x_{j;8} -x_{j+1;7}
& x_{j+2;4} - x_{j+2;6}

 & x_{j+2;8} - x_{j+6;2} \\
 x_{j+3;3} - x_{j+3;8} 

& x_{j+3;6} - x_{j+6;3} 

& x_{j+3;7} - x_{j+6;7}
& x_{j+4;2} - x_{j+4;7} 

&  x_{j+4;8} - x_{j+6;8} \\

 x_{j+5;1} -x_{j+10;1}

 & x_{j+5;5} - x_{j+6;5}

 & x_{j+5;7} - x_{j+10;2}

& x_{j+7;3} - x_{j+7;6}

 & x_{j+7;8} - x_{j+10;3} \\
 x_{j+8;2} - x_{j+8;8} & x_{j+8;6} - x_{j+10;4}

 & x_{j+8;7} - x_{j+10;8}
& x_{j+9;1} - x_{j+9;7} & x_{j+9;8}- x_{j+10;6}\\
 x_{j+11;4} - x_{j+11;5} & x_{j+12;3} - x_{j+12;4} 
& x_{j+13;2} - x_{j+13;3}
& x_{j+14;1} - x_{j+14;2}& -x_{j+15;1}
  \end{array}
\\
  \begin{array}{cc}
 2x_{j+5;5} - x_{j+5;6}-x_{j+5;8}-x_{j+6;4}
& x_{j+6;5} - x_{j+6;7}-x_{j+6;8}-x_{j+10;1} \\
 x_{j+3;7}+x_{j+4;8}+x_{j+5;1} - x_{j+5;5}
& x_{j+5;6}+x_{j+5;8}+x_{j+6;4} - 2x_{j+6;5}
  \end{array}
\\
  \begin{array}{ccc}
 x_{j;6}+x_{j;8} - x_{j+1;5}
& x_{j;7} + x_{j;8}-  x_{j+1;6}
& x_{j;7} + x_{j+2;4} - x_{j+2;5} \\
 x_{j;7} + x_{j+3;3} - x_{j+3;4}
& x_{j;7} +x_{j+4;2} -x_{j+4;3}
& x_{j;7} +x_{j+5;1} - x_{j+5;2} \\

 x_{j+2;8} + x_{j+3;3} - x_{j+3;5}
& x_{j+2;8} + x_{j+4;2} - x_{j+4;4}
& x_{j+2;8} +x_{j+5;1} - x_{j+5;3} \\

 x_{j+3;6} +x_{j+4;2} - x_{j+4;5}
& x_{j+3;6} + x_{j+5;1} - x_{j+5;4}
& x_{j+3;7} + x_{j+4;2} - x_{j+4;6} \\

 x_{j+3;7} +x_{j+4;8} - x_{j+6;4}
& x_{j+3;7} + x_{j+5;1} - x_{j+5;8}
& x_{j+3;7} + x_{j+5;6} - x_{j+6;5} \\
 x_{j+3;7} +x_{j+5;7} - x_{j+6;6}

& x_{j+4;8} + x_{j+5;1} - x_{j+5;6}
& x_{j+4;8} + x_{j+5;8} - x_{j+6;5} \\
 x_{j+5;1} +x_{j+6;4} - x_{j+6;5}
& x_{j+5;1} + x_{j+7;3} - x_{j+7;4}
& x_{j+5;1} + x_{j+8;2} - x_{j+8;3} \\
 x_{j+5;1} +x_{j+9;1} - x_{j+9;2}

& x_{j+5;7} + x_{j+7;3} - x_{j+7;5}
& x_{j+5;7} + x_{j+8;2} - x_{j+8;4} \\
 x_{j+5;7} +x_{j+9;1} - x_{j+9;3}

& x_{j+7;8} + x_{j+8;2} - x_{j+8;5}

& x_{j+7;8} + x_{j+9;1} - x_{j+9;4} \\
 x_{j+8;6} +x_{j+9;1} - x_{j+9;5}

& x_{j+8;7} + x_{j+9;1} - x_{j+9;6}
& x_{j+8;7} + x_{j+9;8} - x_{j+10;5} \\
 x_{j+1;5}  - x_{j+1;7} -x_{j+1;8}
& x_{j+1;6}  - x_{j+1;7} -x_{j+6;1}

& x_{j+2;5}  -x_{j+2;6} - x_{j+6;1} \\
 x_{j+3;4} - x_{j+3;8} -x_{j+6;1}
& x_{j+3;5} - x_{j+3;8} - x_{j+6;2}

& x_{j+4;3}  - x_{j+4;7} - x_{j+6;1} \\
 x_{j+4;4} - x_{j+4;7} -x_{j+6;2}
& x_{j+4;5} - x_{j+4;7} - x_{j+6;3}

& x_{j+4;6}  - x_{j+4;7} - x_{j+6;7} \\

 x_{j+5;2} - x_{j+6;1} -x_{j+10;1}
& x_{j+5;3} - x_{j+6;2} -x_{j+10;1}
& x_{j+5;4} - x_{j+6;3} -x_{j+10;1} \\
 x_{j+5;5} - x_{j+5;6} -x_{j+6;7}

& x_{j+5;5} - x_{j+6;4} -x_{j+10;1}
& x_{j+5;5} - x_{j+5;8} -x_{j+6;8} \\
 x_{j+5;6} - x_{j+6;8} -x_{j+10;1}

& x_{j+5;8} - x_{j+6;7} -x_{j+10;1}
& x_{j+6;4} - x_{j+6;7} -x_{j+6;8} \\

 x_{j+6;6} - x_{j+6;7} -x_{j+10;2}
& x_{j+7;4} - x_{j+7;6} -x_{j+10;1}
& x_{j+7;5} - x_{j+7;6} -x_{j+10;2} \\

 x_{j+8;3} - x_{j+8;8} -x_{j+10;1}
& x_{j+8;4} - x_{j+8;8} -x_{j+10;2}
& x_{j+8;5} - x_{j+8;8} -x_{j+10;3} \\

 x_{j+9;2} - x_{j+9;7} -x_{j+10;1}
& x_{j+9;3} - x_{j+9;7} -x_{j+10;2}
& x_{j+9;4} - x_{j+9;7} -x_{j+10;3} \\
 x_{j+9;5} - x_{j+9;7} -x_{j+10;4}
& x_{j+9;6} - x_{j+9;7} -x_{j+10;8}
& x_{j+10;5} - x_{j+10;6} -x_{j+10;8}
  \end{array}
\\
  \begin{array}{cc}
 x_{j+3;5}+x_{j+4;3}+x_{j+5;1} - x_{j+3;8}-x_{j+4;4}-x_{j+5;2}

& x_{j+3;6}+x_{j+4;4}+x_{j+5;2} - x_{j+4;5}-x_{j+5;3}-x_{j+6;1} \\
 x_{j+3;7}+x_{j+4;4}+x_{j+5;2} -x_{j+4;6}-x_{j+5;3}-x_{j+6;1}
& x_{j+3;7}+x_{j+4;5}+x_{j+5;2} - x_{j+4;6}-x_{j+5;4}-x_{j+6;1} \\
 x_{j+3;7}+x_{j+4;5}+x_{j+5;3} - x_{j+4;6}-x_{j+5;4}-x_{j+6;2}

& x_{j+4;6}+x_{j+4;8}+x_{j+5;2} - x_{j+4;7}-x_{j+5;5}-x_{j+6;1} \\

x_{j+4;6}+x_{j+4;8}+x_{j+5;3} - x_{j+4;7}-x_{j+5;5}-x_{j+6;2}
& x_{j+4;6}+x_{j+4;8}+x_{j+5;4} - x_{j+4;7}-x_{j+5;5}-x_{j+6;3} \\

 x_{j+5;7}+x_{j+6;5}+x_{j+7;3} - x_{j+6;6}-x_{j+6;8}-x_{j+7;4}

& x_{j+5;7}+x_{j+6;5}+x_{j+8;2} - x_{j+6;6}-x_{j+6;8}-x_{j+8;3} \\
 x_{j+5;7}+x_{j+6;5}+x_{j+9;1} - x_{j+6;6}-x_{j+6;8}-x_{j+9;2}

& x_{j+6;6}+x_{j+7;4}+x_{j+8;2} - x_{j+6;7}-x_{j+7;5}-x_{j+8;3} \\
 x_{j+6;6}+x_{j+7;4}+x_{j+9;1} - x_{j+6;7}-x_{j+7;5}-x_{j+9;2}
& x_{j+6;6}+x_{j+8;3}+x_{j+9;1} - x_{j+6;7}-x_{j+8;4}-x_{j+9;2} \\
 x_{j+7;5}+x_{j+8;3}+x_{j+9;1} - x_{j+7;6}-x_{j+8;4}-x_{j+9;2}

& x_{j+7;8}+x_{j+8;4}+x_{j+9;2} - x_{j+8;5}-x_{j+9;3}-x_{j+10;1}
  \end{array}
\end{array}
\]
\[
\begin{array}{cc}
  \begin{array}{ccc}
 x_{j;7}+x_{j+1;5} - x_{j+1;6}-x_{j+1;8}
& x_{j+1;6}+x_{j+2;4} - x_{j+1;7}-x_{j+2;5}
& x_{j+1;6}+x_{j+3;3} - x_{j+1;7}-x_{j+3;4} \\
 x_{j+1;6}+x_{j+4;2} - x_{j+1;7}-x_{j+4;3}
& x_{j+1;6}+x_{j+5;1} - x_{j+1;7}-x_{j+5;2}
& x_{j+2;5}+x_{j+2;6} - x_{j+3;3}-x_{j+3;4} \\

 x_{j+2;5}+x_{j+4;2} -x_{j+2;6}-x_{j+4;3}
& x_{j+2;5}+x_{j+5;1} - x_{j+2;6}-x_{j+5;2}
& x_{j+2;8}+x_{j+3;4} - x_{j+3;5}-x_{j+6;1} \\
 x_{j+2;8}+x_{j+4;3} - x_{j+4;4}-x_{j+6;1}
& x_{j+2;8}+x_{j+5;2} - x_{j+5;3}-x_{j+6;1}
& x_{j+3;4}+x_{j+4;2} - x_{j+3;8}-x_{j+4;3} \\
 x_{j+3;4}+x_{j+5;1} - x_{j+3;8}-x_{j+5;2}
& x_{j+3;5}+x_{j+4;2} - x_{j+3;8}-x_{j+4;4}
& x_{j+3;5}+x_{j+5;1} - x_{j+3;8}-x_{j+5;3} \\

 x_{j+3;6}+x_{j+4;3} - x_{j+4;5}-x_{j+6;1}
& x_{j+3;6}+x_{j+4;4} - x_{j+4;5}-x_{j+6;2}
& x_{j+3;6}+x_{j+5;2} - x_{j+5;4}-x_{j+6;1} \\
 x_{j+3;6}+x_{j+5;3} - x_{j+5;4}-x_{j+6;2}
& x_{j+3;7}+x_{j+4;3} - x_{j+4;6}-x_{j+6;1}
& x_{j+3;7}+x_{j+4;4} - x_{j+4;6}-x_{j+6;2} \\
 x_{j+3;7}+x_{j+4;5} - x_{j+4;6}-x_{j+6;3}


& x_{j+3;7}+x_{j+5;2} - x_{j+5;8}-x_{j+6;1}
& x_{j+3;7}+x_{j+5;3} - x_{j+5;8}-x_{j+6;2} \\
 x_{j+3;7}+x_{j+5;4} - x_{j+5;8}-x_{j+6;3}
& x_{j+3;7}+x_{j+5;5} - x_{j+5;8}-x_{j+6;4}
& x_{j+4;3}+x_{j+5;1} - x_{j+4;7}-x_{j+5;2} \\
 x_{j+4;4}+x_{j+5;1} - x_{j+4;7}-x_{j+5;3}
& x_{j+4;5}+x_{j+5;1} -x_{j+4;7}-x_{j+5;4}

& x_{j+4;6}+x_{j+4;8} - x_{j+4;7}-x_{j+6;4} \\
 x_{j+4;6}+x_{j+5;1} -x_{j+4;7}-x_{j+5;8}
& x_{j+4;6}+x_{j+5;6} - x_{j+4;7}-x_{j+6;5}
& x_{j+4;6}+x_{j+5;7} - x_{j+4;7}-x_{j+6;6} \\
 x_{j+4;8}+x_{j+5;2} - x_{j+5;6}-x_{j+6;1}
& x_{j+4;8}+x_{j+5;3} - x_{j+5;6}-x_{j+6;2}
& x_{j+4;8}+x_{j+5;4} - x_{j+5;6}-x_{j+6;3} \\
 x_{j+4;8}+x_{j+5;5} - x_{j+5;6}-x_{j+6;4}
& x_{j+5;1}+x_{j+5;5} - x_{j+5;6}-x_{j+5;8}

& x_{j+5;2}+x_{j+6;4} - x_{j+6;1}-x_{j+6;5} \\
 x_{j+5;2}+x_{j+7;3} - x_{j+6;1}-x_{j+7;4}
& x_{j+5;2}+x_{j+8;2} - x_{j+6;1}-x_{j+8;3}
& x_{j+5;2}+x_{j+9;1} - x_{j+6;1}-x_{j+9;2} \\
 x_{j+5;3}+x_{j+6;4} - x_{j+6;2}-x_{j+6;5}
& x_{j+5;3}+x_{j+7;3} - x_{j+6;2}-x_{j+7;4}
& x_{j+5;3}+x_{j+8;2} - x_{j+6;2}-x_{j+8;3} \\
 x_{j+5;3}+x_{j+9;1} - x_{j+6;2}-x_{j+9;2}
& x_{j+5;4}+x_{j+6;4} - x_{j+6;3}-x_{j+6;5}
& x_{j+5;4}+x_{j+7;3} - x_{j+6;3}-x_{j+7;4} \\
 x_{j+5;4}+x_{j+8;2} - x_{j+6;3}-x_{j+8;3}
& x_{j+5;4}+x_{j+9;1} -x_{j+6;3}-x_{j+9;2}
& x_{j+5;5}+x_{j+5;7} - x_{j+5;6}-x_{j+6;6} \\

 x_{j+5;5}+x_{j+7;3} - x_{j+6;4}-x_{j+7;4}
& x_{j+5;5}+x_{j+8;2} - x_{j+6;4}-x_{j+8;3}
& x_{j+5;5}+x_{j+9;1} - x_{j+6;4}-x_{j+9;2} \\
 x_{j+5;6}+x_{j+5;8} - x_{j+6;5}-x_{j+10;1}
& x_{j+5;6}+x_{j+6;4} - x_{j+6;5}-x_{j+6;8}
& x_{j+5;6}+x_{j+7;3} - x_{j+6;8}-x_{j+7;4} \\
 x_{j+5;6}+x_{j+8;2} - x_{j+6;8}-x_{j+8;3}
& x_{j+5;6}+x_{j+9;1} - x_{j+6;8}-x_{j+9;2}
& x_{j+5;7}+x_{j+5;8} - x_{j+6;6}-x_{j+10;1} \\
 x_{j+5;7}+x_{j+6;4} - x_{j+6;6}-x_{j+6;8}

& x_{j+5;7}+x_{j+7;4} -x_{j+7;5}-x_{j+10;1}
& x_{j+5;7}+x_{j+8;3} - x_{j+8;4}-x_{j+10;1} \\
 x_{j+5;7}+x_{j+9;2} - x_{j+9;3}-x_{j+10;1}
& x_{j+5;8}+x_{j+6;4} - x_{j+6;5}-x_{j+6;7}
& x_{j+5;8}+x_{j+7;3} -x_{j+6;7}-x_{j+7;4} \\
 x_{j+5;8}+x_{j+8;2} - x_{j+6;7}-x_{j+8;3}
& x_{j+5;8}+x_{j+9;1} - x_{j+6;7}-x_{j+9;2}

& x_{j+6;6}+x_{j+7;3} - x_{j+6;7}-x_{j+7;5} \\
 x_{j+6;6}+x_{j+8;2} - x_{j+6;7}-x_{j+8;4}
& x_{j+6;6}+x_{j+9;1} - x_{j+6;7}-x_{j+9;3}
& x_{j+7;4}+x_{j+8;2} - x_{j+7;6}-x_{j+8;3} \\
 x_{j+7;4}+x_{j+9;1} - x_{j+7;6}-x_{j+9;2}
& x_{j+7;5}+x_{j+8;2} - x_{j+7;6}-x_{j+8;4}
& x_{j+7;5}+x_{j+9;1} - x_{j+7;6}-x_{j+9;3} \\

 x_{j+7;8}+x_{j+8;3} - x_{j+8;5}-x_{j+10;1}
& x_{j+7;8}+x_{j+8;4} - x_{j+8;5}-x_{j+10;2}
& x_{j+7;8}+x_{j+9;2} - x_{j+9;4}-x_{j+10;1} \\
 x_{j+7;8}+x_{j+9;3} - x_{j+9;4}-x_{j+10;2}
& x_{j+8;3}+x_{j+9;1} - x_{j+8;8}-x_{j+9;2}
& x_{j+8;4}+x_{j+9;1} - x_{j+8;8}-x_{j+9;3} \\
 x_{j+8;5}+x_{j+9;1} - x_{j+8;8}-x_{j+9;4}
& x_{j+8;6}+x_{j+9;2} - x_{j+9;5}-x_{j+10;1}
& x_{j+8;6}+x_{j+9;3} - x_{j+9;5}-x_{j+10;2} \\
 x_{j+8;6}+x_{j+9;4} -x_{j+9;5}-x_{j+10;3}

& x_{j+8;7}+x_{j+9;2} -x_{j+9;6}-x_{j+10;1}
& x_{j+8;7}+x_{j+9;3} - x_{j+9;6}-x_{j+10;2} \\
 x_{j+8;7}+x_{j+9;4} - x_{j+9;6}-x_{j+10;3}
& x_{j+8;7}+x_{j+9;5} - x_{j+9;6}-x_{j+10;4}
& x_{j+9;6}+x_{j+9;8} -x_{j+9;7}-x_{j+10;5}
  \end{array}
\\
  \begin{array}{cc}
 x_{j+2;8}+x_{j+3;4}+x_{j+4;2} - x_{j+3;5}-x_{j+4;3}
& x_{j+2;8}+x_{j+3;4}+x_{j+5;1} -x_{j+3;5}-x_{j+5;2} \\
 x_{j+2;8}+x_{j+4;3}+x_{j+5;1} -x_{j+4;4}-x_{j+5;2}

& x_{j+3;6}+x_{j+4;3}+x_{j+5;1} -x_{j+4;5}-x_{j+5;2} \\
 x_{j+3;6}+x_{j+4;4}+x_{j+5;1} - x_{j+4;5}-x_{j+5;3}
& x_{j+3;7}+x_{j+4;3}+x_{j+5;1} - x_{j+4;6}-x_{j+5;2} \\
 x_{j+3;7}+x_{j+4;4}+x_{j+5;1} - x_{j+4;6}-x_{j+5;3}
& x_{j+3;7}+x_{j+4;5}+x_{j+5;1} - x_{j+4;6}-x_{j+5;4} \\
 x_{j+3;7}+x_{j+4;8}+x_{j+5;2} - x_{j+5;5}-x_{j+6;1}

& x_{j+3;7}+x_{j+4;8}+x_{j+5;3} - x_{j+5;5}-x_{j+6;2} \\
 x_{j+3;7}+x_{j+4;8}+x_{j+5;4} - x_{j+5;5}-x_{j+6;3}
& x_{j+4;6}+x_{j+4;8}+x_{j+5;1} - x_{j+4;7}-x_{j+5;5} \\

 x_{j+5;6}+x_{j+5;8}+x_{j+7;3} - x_{j+6;5}-x_{j+7;4}
& x_{j+5;6}+x_{j+5;8}+x_{j+8;2} - x_{j+6;5}-x_{j+8;3} \\
 x_{j+5;6}+x_{j+5;8}+x_{j+9;1} - x_{j+6;5}-x_{j+9;2}
& x_{j+5;7}+x_{j+5;8}+x_{j+6;4} - x_{j+6;5}-x_{j+6;6} \\
 x_{j+5;7}+x_{j+5;8}+x_{j+7;3} - x_{j+6;6}-x_{j+7;4}
& x_{j+5;7}+x_{j+5;8}+x_{j+8;2} - x_{j+6;6}-x_{j+8;3} \\
 x_{j+5;7}+x_{j+5;8}+x_{j+9;1} - x_{j+6;6}-x_{j+9;2}

& x_{j+5;7}+x_{j+7;4}+x_{j+8;2} -x_{j+7;5}-x_{j+8;3} \\
 x_{j+5;7}+x_{j+7;4}+x_{j+9;1} - x_{j+7;5}-x_{j+9;2}
& x_{j+5;7}+x_{j+8;3}+x_{j+9;1} - x_{j+8;4}-x_{j+9;2} \\

 x_{j+7;8}+x_{j+8;3}+x_{j+9;1} -x_{j+8;5}-x_{j+9;2}

& x_{j+7;8}+x_{j+8;4}+x_{j+9;1} - x_{j+8;5}-x_{j+9;3} \\
 x_{j+3;5}+x_{j+4;3} - x_{j+3;8}-x_{j+4;4}-x_{j+6;1}
& x_{j+3;5}+x_{j+5;2} - x_{j+3;8}-x_{j+5;3}-x_{j+6;1} \\

 x_{j+4;4}+x_{j+5;2} - x_{j+4;7}-x_{j+5;3}-x_{j+6;1}
& x_{j+4;5}+x_{j+5;2} -x_{j+4;7}-x_{j+5;4}-x_{j+6;1} \\
 x_{j+4;5}+x_{j+5;3} - x_{j+4;7}-x_{j+5;4}-x_{j+6;2}

& x_{j+4;6}+x_{j+5;2} -x_{j+4;7}-x_{j+5;8}-x_{j+6;1} \\
 x_{j+4;6}+x_{j+5;3} - x_{j+4;7}-x_{j+5;8}-x_{j+6;2}
& x_{j+4;6}+x_{j+5;4} - x_{j+4;7}-x_{j+5;8}-x_{j+6;3} \\
 x_{j+4;6}+x_{j+5;5} - x_{j+4;7}-x_{j+5;8}-x_{j+6;4}

& x_{j+5;2}+x_{j+5;5} - x_{j+5;6}-x_{j+5;8}-x_{j+6;1} \\
 x_{j+5;3}+x_{j+5;5} -x_{j+5;6}-x_{j+5;8}-x_{j+6;2}
& x_{j+5;4}+x_{j+5;5} - x_{j+5;6}-x_{j+5;8}-x_{j+6;3} \\

 x_{j+5;7}+x_{j+6;5} - x_{j+6;6}-x_{j+6;8}-x_{j+10;1}
& x_{j+6;5}+x_{j+7;3} - x_{j+6;7}-x_{j+6;8}-x_{j+7;4} \\
 x_{j+6;5}+x_{j+8;2} - x_{j+6;7}-x_{j+6;8}-x_{j+8;3}
& x_{j+6;5}+x_{j+9;1} - x_{j+6;7}-x_{j+6;8}-x_{j+9;2} \\

 x_{j+6;6}+x_{j+7;4} - x_{j+6;7}-x_{j+7;5}-x_{j+10;1}
& x_{j+6;6}+x_{j+8;3} - x_{j+6;7}-x_{j+8;4}-x_{j+10;1} \\
 x_{j+6;6}+x_{j+9;2} - x_{j+6;7}-x_{j+9;3}-x_{j+10;1}
& x_{j+7;5}+x_{j+8;3} - x_{j+7;6}-x_{j+8;4}-x_{j+10;1} \\
 x_{j+7;5}+x_{j+9;2} - x_{j+7;6}-x_{j+9;3}-x_{j+10;1}

& x_{j+8;4}+x_{j+9;2} - x_{j+8;8}-x_{j+9;3}-x_{j+10;1} \\
 x_{j+8;5}+x_{j+9;2} - x_{j+8;8}-x_{j+9;4}-x_{j+10;1}
& x_{j+8;5}+x_{j+9;3} -x_{j+8;8}-x_{j+9;4}-x_{j+10;2}
  \end{array}
\end{array}
\]
\end{thm1}

\noindent
(i) By the form of $\Sigma_{\io}$, coefficients of $x_{1;1}$, $x_{1;2}$, 
$x_{1;3}$, $x_{1;4}$, $x_{1;5}$, $x_{1;6}$, $x_{1;7}$, $x_{1;8}$
are positive. This shows that ``positivity 
assumption'' is satisfied.
\vskip5pt

\noindent
$(ii)$
We assume $x_{j;1} \geq 0$ for any $j \geq 1$.
We show that $x_{j;2}$, $x_{j;3}$, $x_{j;4}$,
$x_{j;5}$, $x_{j;6}$, $x_{1;7}$, $x_{1;8} \geq 0$ for any $j \geq 1$.
By the form of $\Sigma_{\io}$, we have
$$x_{j;2} \geq x_{j+1;1},\;\;
x_{j;3} \geq x_{j+1;2},\;\;
x_{j;4} \geq x_{j+1;3},\;\;
x_{j;5} \geq x_{j+1;4},\;\;
x_{j;7} \geq x_{j+6;1},\;\;
x_{j;8} \geq x_{j+1;7},\;\;
x_{j;6} \geq x_{j+1;8}.$$
This shows $x_{j;2} \geq 0$ for $ j\geq 1$ since $x_{j+1;1} \geq 0$ and
similarly, we have $x_{j;3} \geq 0$, $x_{j;4} \geq 0$,
$x_{j;5} \geq 0$, $x_{j;7} \geq 0$, $x_{j;8} \geq 0$,
$x_{j;6} \geq 0$
for $j \geq 1$
since $x_{j;3} \geq x_{j+1;2}$,
$x_{j;4} \geq x_{j+1;3}$, $x_{j;5} \geq x_{j+1;4}$,
$x_{j;7} \geq x_{j+1;6}$, $x_{j;8} \geq x_{j+1;7}$,
$x_{j;6} \geq x_{j+1;8}$ respectively.
\vskip5pt















\noindent
$(iii)$
We determine when $x_{j;i} \equiv 0$.
We have
$$0 \geq x_{j+15;1} \geq 0.$$
This shows $x_{m;1} \equiv 0$ for $m \geq 16$.
Similarly, we have
$x_{m;2} \equiv 0$, 
$x_{m;3} \equiv 0$,
$x_{m;4} \equiv 0$,
$x_{m;5} \equiv 0$,
$x_{m;6} \equiv 0$,
$x_{m;7} \equiv 0$,
$x_{m;8} \equiv 0$
for $m \geq 16$
since $x_{j+14;1} \geq x_{j+14;2}$,
$x_{j+13;2} \geq x_{j+13;3}$,
$x_{j+12;3} \geq x_{j+12;4}$,
$x_{j+11;4} \geq x_{j+11;5}$,
$x_{j+7;3} \geq x_{j+7;6}$,
$x_{j+9;1} \geq x_{j+9;7}$,
$x_{j+8;2} \geq x_{j+8;8}$
respectively. In particular, the parameter $j$ in $\Xi_{\io}$ runs
$1 \leq j \leq 15$.
\begin{rem1}
For all simple case,
the number of $x_{j;i}$'s such that $x_{j;i} \not\equiv 0$
are equal to the length of longest element
of Weyl group.
\label{number}
\end{rem1}




\begin{thebibliography}{99}
\def\CMP{\sl Commun.Math.Phys.}
\def\IJMP{\sl Int.J.Mod.Phys.}
\def\Duke{\sl Duke Math.J.}


%



\bibitem[H]{H} A. Hoshino, Polyhedral Realizations of
             Crystal Bases for Modified Quantum Algebras of Arbitrary 
             Rank 2 Cases, math.QA/0403192.

\bibitem[HN]{HN} A. Hoshino and T. Nakashima, Polyhedral Realizations of
             Crystal Bases for Modified Quantum Algebras of Type $A$, 
             to appear Comm in Algebra.
 
\bibitem[Kac]{Kac} V. Kac, Infinite dimensional Lie algebras, 
	      3rd edition, Cambridge Univ. Press, 1990.


\bibitem[K1]{K1} M. Kashiwara,
 On crystal bases of the $q$-analogue of universal enveloping algebras,
	Duke Math. J., {\bf 63} (1991), 465--516.


\bibitem[K2]{K2} M. Kashiwara,
      Crystal base and Littelmann's refined Demazure character formula,
	     Duke Math. J., {\bf 71} (1993), 839--858.











\bibitem[Li]{Li} P. Littelmann, Corns, crystals, and patterns,
Transform. Groups, {\bf 3} (1998), no. 2, 145-179.

\bibitem[L1]{L1} G. Lusztig, Canonical bases in tensor products,
Proc. Nat. Acad. Sci. U.S.A., {\bf 89} (1992), 8177-8179.

\bibitem[L2]{L2} G. Lusztig, Introduction to quantum groups,
Birkh\"auser, Boston, 1993.



\bibitem[N]{N1} T. Nakashima,  Polyhedral Realizations of 
       Crystal Bases for Integrable Highest Weight Modules, 
           J. Algebra, {\bf 219} (1999), 571--597.

\bibitem[NZ]{NZ} T. Nakashima and A. Zelevinsky,  Polyhedral Realization of 
       Crystal Bases for Quantized Kac-Moody Algebras, 
           Adv. in Math., {\bf 131} No.1 (1997), 253--278.

\end{thebibliography}
\end{document}